# THE SPECTRUM OF KERNEL RANDOM MATRICES

By Noureddine El Karoui[1]

*University of California, Berkeley*

We place ourselves in the setting of high-dimensional statistical inference where the number of variables $p$ in a dataset of interest is of the same order of magnitude as the number of observations $n$.

We consider the spectrum of certain kernel random matrices, in particular $n \times n$ matrices whose $(i,j)$th entry is $f(X_i' X_j/p)$ or $f(\|X_i - X_j\|^2/p)$ where $p$ is the dimension of the data, and $X_i$ are independent data vectors. Here $f$ is assumed to be a locally smooth function.

The study is motivated by questions arising in statistics and computer science where these matrices are used to perform, among other things, nonlinear versions of principal component analysis. Surprisingly, we show that in high-dimensions, and for the models we analyze, the problem becomes essentially linear—which is at odds with heuristics sometimes used to justify the usage of these methods. The analysis also highlights certain peculiarities of models widely studied in random matrix theory and raises some questions about their relevance as tools to model high-dimensional data encountered in practice.

**1. Introduction.** Recent years has seen newfound theoretical interest in the properties of large-dimensional sample covariance matrices. With the increase in the size and dimensionality of datasets to be analyzed, questions have been raised about the practical relevance of information derived from classical asymptotic results concerning spectral properties of sample covariance matrices. To address these concerns, one line of analysis has been the consideration of asymptotics where both the sample size $n$ and the number of variables $p$ in the dataset go to infinity jointly while assuming, for instance, that $p/n$ had a limit.

Received December 2007; revised August 2008.
[1]Supported by NSF Grant DMS-06-05169.

This type of questions concerning the spectral properties of large-dimensional matrices have been and are being addressed in variety of fields, from physics to various areas of mathematics. While the topic is classical, with the seminal contribution [43] dating back from the 1950s, there has been renewed and vigorous interest in the study of large-dimensional random matrices in the last decade or so. This has led to new insights and the appearance of new "canonical" distributions [38], new tools (see [41]) and, in Statistics, a sense that one needs to exert caution with familiar techniques of multivariate analysis when the dimension of the data gets large and the sample size is of the same order of magnitude as that dimension.

So far in Statistics, this line of work has been concerned mostly with the properties of sample covariance matrices. In a seminal paper, Marčenko and Pastur [30] showed a result that, from a statistical standpoint, may be interpreted as saying, roughly, that asymptotically, the histogram the eigenvalues of a sample (i.e., random) covariance matrix is (asymptotically) a *deterministic* nonlinear deformation of the histogram of the eigenvalues of the population covariance matrix. Remarkably, they managed to characterize this deformation for fairly general population covariances. Their result was shown in great generality and introduced new tools to the field including one that has become ubiquitous, the Stieltjes transform of a distribution. In its best-known form, their result says that when the population covariance is identity, and hence all the population eigenvalues are equal to 1, in the limit the sample eigenvalues are split, and if $p \leq n$, they are spread between $[(1-\sqrt{p/n})^2, (1+\sqrt{p/n})^2]$ according to a fully explicit density known now as the density of the Marčenko–Pastur law. Their result was later re-discovered independently in [42] (under slightly weaker conditions) and generalized to the case of nondiagonal covariance matrices in [37] under some particular distributional assumptions which we discuss later in the paper.

On the other hand, recent developments have been concerned with fine properties of the largest eigenvalue of random matrices which became amenable to analysis after mathematical breakthroughs which happened in the 1990s (see [38, 39] and [40]). Classical statistical work on joint distribution of eigenvalues of sample covariance matrices (see [1] for a good reference) then became usable for analysis in high-dimensions. In particular, in the case of gaussian distributions, with Id covariance, it was shown in [27] and [16] that the largest eigenvalue of the sample covariance matrix is Tracy–Widom distributed. More recent progress [17] managed to carry out the analysis for essentially general population covariance. On the other hand, models for which the population covariance has a few separated eigenvalues have also been of interest (see, for instance, [8] and [31]). Beside the particulars of the different type of fluctuations that can be encountered (Tracy–Widom, Gaussian or other), researchers have been able to precisely localize these largest eigenvalues. One interesting aspect of those results is the fact that in the high-dimensional setting of interest to us, the largest eigenvalues are always positively biased, with the bias being sometime large. (We also note that in the case of i.i.d. data—which naturally is less interesting in statistics—results on the localization of the largest eigenvalue have been available for quite some time now, after the works [21] and [45], to cite a few.) This is naturally in sharp contrast to classical results of multivariate analysis which show $\sqrt{n}$-consistency of all sample eigenvalues; though the possibility of bias is a simple consequence of Jensen's inequality.

On the other hand, there has been much less theoretical work on kernel random matrices. By this term, we mean matrices with $(i,j)$ entry of the form,

$$M_{i,j} = k(X_i, X_j),$$

where $M$ is an $n \times n$ matrix, $X_i$ is a $p$-dimensional data vector and $k$ is a function of two variables, often called a kernel, that may depend on $n$. Common choices of kernels include, for instance, $k(X_i, X_j) = f(\|X_i - X_j\|^2/t)$, where $f$ is a function and $t$ is a scalar, or $k(X_i, X_j) = f(X_i'X_j/t)$. For the function $f$, common choices include $f(x) = \exp(-x)$, $f(x) = \exp(-x^a)$, for a certain scalar $a$, $f(x) = (1+x)^a$, or $f(x) = \tanh(a+bx)$, where $b$ is a scalar. We refer the reader to [33], Chapter 4, or [44] for more examples.

In particular, we are not aware of any work in the setting of high-dimensional data analysis, where $p$ grows with $n$. However, given the practical success and flexibility of these methods (we refer to [35] for an introduction), it is natural to try to investigate theoretically their properties. Further, as illustrated in the data analytic part of [44], an $n/p$ boundedness assumption is not unrealistic as far as applications of kernel methods are concerned. One aim of the present paper is to shed some theoretical light on the properties of these kernel random matrices and to do so in relatively wide generality. We note that the choice of renormalization that we make below is motivated in part by the arguments of [44] and their practical choices of kernels for data of varying dimensions.

Existing theory on kernel random matrices (see, for instance, the interesting [28]) for fixed-dimensional input data predicts that the eigenvalues of kernel random matrices behave—at least for the largest ones—like the eigenvalues of the corresponding operator on $L^2(dP)$, if the data is i.i.d. with probability distribution $P$. To be more precise, if $X_i$ is a sequence of i.i.d. random variables with distribution $P$, under regularity conditions on the kernel $k(x,y)$, it was shown in [28] that, for any index $l$, the $l$th largest eigenvalue of the kernel matrix $M$ with entries,

$$M_{i,j} = \frac{1}{n} k(X_i, X_j),$$

converges to the $l$th largest eigenvalue of the operator $K$ defined as

$$Kf(x) = \int k(x,y)f(y)\,dP(y).$$

These insights have also been derived through more heuristic but nonetheless enlightening arguments in, for instance, [44]. Further, more precise fluctuation results are also given in [28]. We also note interesting work on Laplacian eigenmaps (see, e.g., [9]) where, among other things, results have been obtained showing convergence of eigenvalues and eigenvectors of certain Laplacian random matrices (which are quite closely connected to kernel random matrices) computed from data sampled from a manifold, to corresponding quantities for the Laplace–Beltrami operator on the manifold.

These results are in turn used in the literature to explain the behavior of nonlinear versions of standard procedures of multivariate statistics, such as Principal Component Analysis (PCA), Canonical Correlation Analysis (CCA) or Independent Component Analysis (CCA). We refer the reader to [36] for an introduction to kernel-PCA, and to [2] for an introduction to kernel-CCA and kernel-ICA. At the heart of these techniques are the spectral properties of kernel random matrices. Because these techniques are used in bioinformatics, a field where large datasets are common and becoming the norm, it is natural to ask what can be said about these spectral properties for high-dimensional data.

We show that for the models we analyze (ICA-type models and generalizations that go beyond the linear setting of ICA), kernel random matrices essentially behave like sample covariance matrices and hence their eigenvalues suffer from the same bias problems that affect sample covariance matrices in high-dimensions. In particular, if one were to try to apply the heuristics of [44] which were developed for low-dimensional problems, to the high-dimensional case, the predictions would be quite wildly wrong. (A simple example is provided by the Gaussian kernel with i.i.d. Gaussian data where the computations can be done completely explicitly as explained in [44].) We also note that the scaling we use is different from the one used in low dimensions, where the matrices are scaled by $1/n$. This is because the high-dimensional problem would be completely degenerate if we used this normalization in our setting. However, our results still give information about the problem when it is scaled by $1/n$.

From a random matrix point of view, our study is connected to the study of Euclidean random matrices and distance matrices which is of some interest in, for instance, Physics. We refer to [11] and [12] for work in this direction in the low (or fixed) dimensional setting. We also note that at the level of generality we place ourselves in, the random matrices we study do not seem to be amenable to study through the classical tools of random matrix theory.

Hence, beside their obvious statistical interest, they are also interesting on purely mathematical grounds.

We now turn to the gist of our paper, which will show that high-dimensional kernel random matrices behave spectrally essentially like matrices closely connected to sample covariance matrices. We will get two types of results: in Theorems 2.1 and 2.2, we get a strong approximation result (in operator norm) for standard models (ICA-like) studied in random matrix theory. In Theorems 2.3 and 2.4, we characterize the limiting spectral distribution of our kernel random matrices, for a wider class of data distributions. In Section 2, we also state clearly the consequences of our theorems and review the relevant theory of high-dimensional sample covariance matrices. From a technical standpoint, we adopt a point of view centered on the concentration of measure phenomenon, as exposed for instance in [29] as it provides a unified way to treat the two types of results we are interested in. Finally, we discuss in our (self-contained) conclusion (Section 3), the consequences of our results and in particular some possible limitations of "standard" random matrix models as tools to model data encountered in practice focusing on geometric properties of datasets drawn according to those models. As explained in more details there, vectors drawn according to these standard random matrix models essentially live close to spheres and are almost orthogonal to one another, a property that may or may not be present in datasets to be analyzed and can be seen as a key to many classical and less classical random matrix results (see also [19]).

**2. Spectrum of kernel random matrices.** Kernel random matrices do not seem to be amenable to analysis through the usual tools of random matrix theory. In particular, for general $f$, it seems difficult to carry out either a method of moments proof, or a Stieltjes transform proof, or a proof that relies on knowing the density of the eigenvalues of the matrix.

Hence, we take an indirect approach. Our strategy is to find approximations of the kernel random matrix that have two properties. First, the approximation matrix is analyzable or has already been analyzed in random matrix theory. Second, the quality of the approximation is good enough that spectral properties of the approximating matrix can be shown to carry over to the kernel matrix.

The strategy in the first two theorems is to derive an operator norm "consistent" approximation of our kernel matrix. In other words, if we call $M$ our kernel matrix, we will find $K$ such that $\|M - K\|_2 \to 0$, as $n$ and $p$ tend to $\infty$. Note that both $M$ and $K$ are real symmetric (and hence Hermitian) here. We explain after the statement of Theorem 2.1 why operator norm consistency is a desirable property. But let us say that in a nutshell, it implies consistency for each individual eigenvalue as well as eigenspaces corresponding to separated eigenvalues.

For the second set of theorems (Theorems 2.3 and 2.4), we will relax the distributional assumptions made on the data, but, at the expense of the precision of the results we will obtain, we will characterize the limiting spectral distribution of our kernel random matrices.

Our theorems below show that kernel random matrices can be well approximated by matrices that are closely connected to large-dimensional covariance matrices. The spectral properties of those matrices have been the subject of a significant amount of work in recent and less recent years, and hence this knowledge, or at least part of it, can be transferred to kernel random matrices. In particular, we refer the reader to [4, 5, 8, 17, 19, 21, 27, 30, 31, 37, 42] and [45] for some of the most statistically relevant results in this area. We review some of them now.

2.1. *Some results on large-dimensional sample covariance matrices.* Since our main theorems are approximating theorems, we first wish to state some of the properties of the objects we will use to approximate kernel random matrices. In what follows, we consider an $n \times p$ data matrix, with, say $p/n$ having a finite nonzero limit. Most of the results that have been obtained are of two types: either they are so-called "bulk" results and concern essentially the spectral distribution (or loosely speaking the histogram of eigenvalues) of the random matrices of interest; or they concern the localization and fluctuation behavior of extreme eigenvalues of these random matrices.

2.1.1. *Spectral distribution results.* An object of interest in random matrix theory is the spectral distribution of random matrices. Let us call $l_i$ the decreasingly ordered eigenvalues of our random matrix, and let us assume we are working with an $n \times n$ matrix, $M_n$. The empirical spectral distribution of a $n \times n$ matrix is the probability measure which puts mass $1/n$ at each of its eigenvalues. In other words, if we call $F_n$ this probability measure, we have

$$dF_n(x) = \frac{1}{n} \sum_{i=1}^{n} \delta_{l_i}(x).$$

Note that the histogram of eigenvalues represents an integrated version of this measure.

For random matrices, this measure $F_n$ is naturally a random measure. A key result in the area of covariance matrices is that if we observe i.i.d. data vectors $X_i$, with $X_i = \Sigma_p^{1/2} Y_i$, where $\Sigma_p$ is a positive semi-definite matrix and $Y_i$ is a vector with i.i.d entries, under weak moment conditions on $Y_i$ and assuming that the spectral distribution of $\Sigma_p$ has a limit (in the sense of weak convergence of distributions), $F_n$ converges to a *nonrandom* measure which we call $F$.

We call the models $X_i = \Sigma_p^{1/2} Y_i$ the "standard" models of random matrix theory because most results have been derived under these assumptions. In particular, various results [5, 6, 21] show, among many other things, that when the entries of the vector $Y$ have 4 (absolute) moments, the largest eigenvalues of the sample covariance matrix $X'X/n$, where $X_i$ now occupies the $i$th row of the $n \times p$ matrix $X$, stay close to the endpoint of the support of $F$.

A natural question is therefore to try to characterize $F$. Except in particular situations, it is difficult to do so explicitly. However, it is possible to characterize a certain transformation of $F$. The tool of choice in this context is the *Stieltjes transform* of a distribution. It is a function defined on $\mathbb{C}^+$ by the formula, if we call $\text{St}_F$ the Stieltjes transform of $F$,

$$\text{St}_F(z) = \int \frac{dF(\lambda)}{\lambda - z}, \qquad \text{Im}[z] > 0.$$

In particular for empirical spectral distributions, we see that, if $F_n$ is the spectral distribution of the matrix $M_n$,

$$\text{St}_{F_n}(z) = \frac{1}{n} \sum_{i=1}^n \frac{1}{l_i - z} = \frac{1}{n} \text{trace}((M_n - z\,\text{Id})^{-1}).$$

The importance of the Stieltjes transform in the context of random matrix theory stems from two facts: on the one hand, it is connected fairly explicitly to the matrices that are being analyzed; on the other hand, pointwise convergence of Stieltjes transform implies weak convergence of distributions, if a certain mass preservation condition is satisfied. This is how a number of bulk results are therefore proved. For a clear and self-contained introduction to the connection between Stieltjes transforms and weak convergence of probability measures, we refer the reader to [22].

The result of [30], later generalized by [37] for standard random matrix models with nondiagonal covariance, and more recently by [19], away from those standard models, is a functional characterization of the limit $F$. If one calls $w_n(z)$ the Stieltjes transform of the empirical spectral distribution of $XX'/n$, $w_n(z)$ converges pointwise (and almost surely after [37]) to a *nonrandom* $w(z)$ which, as a function, is a Stieltjes transform. Moreover, $w$, the Stieltjes transform of $F$, satisfies the equation, if $p/n \to \rho$, $\rho > 0$,

$$-\frac{1}{w(z)} = z - \rho \int \frac{\lambda\,dH(\lambda)}{1 + \lambda w},$$

where $H$ is the limiting spectral distribution of $\Sigma_p$, assuming that such a distribution exists. We note that [37] proved the result under a second moment condition on the entries of $Y_i$.

From this result, [30] derived that in the case where $\Sigma_p = \mathrm{Id}$, and hence $dH = \delta_1$, the empirical spectral distribution has a limit whose density is, if $\rho \leq 1$,

$$f_\rho(x) = \frac{1}{2\pi\rho} \frac{\sqrt{(b-x)(x-a)}}{x},$$

where $a = (1 - \rho^{1/2})^2$ and $b = (1 + \rho^{1/2})^2$. The difference between the population spectral distribution (a point mass at 1, of mass 1) and the limit of the empirical spectral distribution is quite striking.

2.1.2. *Largest eigenvalues results.* Another line of work has been focused on the behavior of extreme eigenvalues of sample covariance matrices. In particular, [21] showed, under some moment conditions, that when $\Sigma_p = \mathrm{Id}_p$, $l_1(X'X/n) \to (1 + \sqrt{p/n})^2$ almost surely. In other words, the largest eigenvalue stays close to the endpoint of the limiting spectral distribution of $X'X/n$. This result was later generalized in [45], and was shown to be true under the assumption of finite fourth moment only, for data with mean 0. In recent years, fluctuation results have been obtained for this largest eigenvalue which is of practical interest in Principal Components Analysis (PCA). Under Gaussian assumptions, [16] and [27] (see also [20] and [26]) showed that the fluctuations of the largest eigenvalue are Tracy–Widom distributed. For the general covariance case, similar results, as well as localization information, were recently obtained in [17]. We note that the localization information (i.e., a formula) that was discovered in this latter paper was shown to hold for a wide variety of standard random matrix models through appeal to [5]. We refer the interested reader to Fact 2 in [17] for more information. Interesting work has also been done on so-called "spiked" models where a few population eigenvalues are separated from the bulk of them. In particular, in the case where all population eigenvalues are equal, except for one that is significantly larger (see [7] for the discovery of an interesting phase transition), [31] showed, in the Gaussian case, inconsistency of the largest sample eigenvalue, as well as the fact that the angle between the population and sample principal eigenvectors is bounded away from 0. Paul [31] also obtained fluctuation information about the largest empirical eigenvalue. Finally, we note that the same inconsistency of eigenvalue result was also obtained in [8], beyond the Gaussian case.

2.1.3. *Notation.* Let us now define some notation and add some clarifications.

We denote by $A'$ the transpose of $A$. The matrices we will be working with all have real entries. We remind the reader that if $A$ and $B$ are two rectangular matrices, $AB$ and $BA$ have the same eigenvalues, except for possibly, a certain number of zeros. We will make repeated use of this fact, for example, for matrices like $X'X$ and $XX'$. In the case where $A$ and $B$ are both square, $AB$ and $BA$ have exactly the same eigenvalues.

We will also need various norms on matrices. We will use the so-called operator norm, which we denote by $|\!|\!|A|\!|\!|_2$ which corresponds to the largest singular value of $A$, that is, $\max_i \sqrt{l_i(A'A)}$. We occasionally denote the largest singular value of $A$ by $\sigma_1(A)$. Clearly, for positive semi-definite matrices, the largest singular value is equal to the largest eigenvalue. Finally, we will sometimes need to use the Frobenius (or Hilbert–Schmidt) norm of a matrix $A$. We denote it by $\|A\|_F$. By definition, it is simply, because we are working with matrices with real entries,

$$\|A\|_F^2 = \sum_{i,j} A_{i,j}^2.$$

Further, we use $\circ$ to denote the Hadamard (i.e., entrywise) product of two matrices. We denote by $\mu_m$ the $m$th moment of a random variable. Note that by a slight abuse of notation, we might also use the same notation to refer to the $m$th absolute moment (i.e., $E|X|^m$) of a random variable, but if there is any ambiguity, we will naturally make precise which definition we are using.

Finally, in the discussion of standard random matrix models that follows, there will be arrays of random variables and a.s. convergence. We work with random variables defined on a common probability space. To each $\omega$ corresponds an infinite-dimensional array of numbers. Unless otherwise noted, the $n \times p$ matrices we will use in what follows are the "upper-left" corner of this array.

We now turn to the study of kernel random matrices. We will show that we can approximate them by matrices that are closely connected to sample covariance matrices in high-dimensions and, therefore, that a number of the results we just reviewed also apply to them.

2.2. *Inner-product kernel matrices:* $f(X_i'X_j/p)$.

THEOREM 2.1 (Spectrum of inner product kernel random matrices). *Let us assume that we observe $n$ i.i.d. random vectors, $X_i$ in $\mathbb{R}^p$. Let us consider the kernel matrix $M$ with entries*

$$M_{i,j} = f\left(\frac{X_i'X_j}{p}\right).$$

*We assume that:*

(a) $n \asymp p$, *that is, $n/p$ and $p/n$ remain bounded as $p \to \infty$.*
(b) $\Sigma_p$ *is a positive semi-definite $p \times p$ matrix, and $|\!|\!|\Sigma_p|\!|\!|_2 = \sigma_1(\Sigma_p)$ remains bounded in $p$, that is, there exists $K > 0$, such that $\sigma_1(\Sigma_p) \leq K$, for all $p$.*

(c) $\Sigma_p/p$ has a finite limit, that is, there exists $l \in \mathbb{R}$ such that $\lim_{p \to \infty} \operatorname{trace}(\Sigma_p)/p = l$.
(d) $X_i = \Sigma_p^{1/2} Y_i$.
(e) *The entries of $Y_i$, a p-dimensional random vector, are i.i.d. Also, denoting by $Y_i(k)$ the kth entry of $Y_i$, we assume that $\mathbf{E}(Y_i(k)) = 0$, $\operatorname{var}(Y_i(k)) = 1$ and $\mathbf{E}(|Y_i(k)|^{4+\varepsilon}) < \infty$ for some $\varepsilon > 0$. (We say that $Y_i$ has $4 + \varepsilon$ absolute moments.)*
(f) *$f$ is a $C^1$ function in a neighborhood of $l = \lim_{p \to \infty} \operatorname{trace}(\Sigma_p)/p$ and a $C^3$ function in a neighborhood of $0$.*

*Under these assumptions, the kernel matrix $M$ can (in probability) be approximated consistently in operator norm, when $p$ and $n$ tend to $\infty$, by the matrix $K$, where*

$$K = \left( f(0) + f''(0) \frac{\operatorname{trace}(\Sigma_p^2)}{2p^2} \right) 11' + f'(0) \frac{XX'}{p} + v_p \operatorname{Id}_n,$$

*where*

$$v_p = f\left( \frac{\operatorname{trace}(\Sigma_p)}{p} \right) - f(0) - f'(0) \frac{\operatorname{trace}(\Sigma_p)}{p}.$$

*In other words,*

$$\|\|M - K\|\|_2 \to 0 \qquad \text{in probability, when } p \to \infty.$$

The advantages of obtaining an operator norm consistent estimator are many. We list some here:

- Asymptotically, $M$ and $K$ have the same $j$-largest eigenvalue, for any $j$; this is simply because for symmetric matrices, if $l_j$ is the $j$th largest eigenvalue of a matrix, Weyl's inequality (see, e.g., Corollary III.2.6 in [10]) implies that

$$|l_j(M) - l_j(K)| \leq \|\|M - K\|\|_2.$$

Hence our result implies that $|l_j(M) - l_j(K)| \to 0$ in probability as $p$ and $n$ go to infinity.
- The limiting spectral distributions of $M$ and $K$ (if they exist) are the same. This is a consequence of Lemma 2.1, page 31 below. So in particular, when $K$ has a limiting spectral distribution (in the sense of weak convergence of probability measures), the empirical spectral distribution of $M$ converges to that distribution (in the sense of weak convergence of distributions) in probability.
- We have subspace consistency for eigenspaces corresponding to separated eigenvalues. (For a proof, we refer to [18], Corollary 3.) So, when $K$ has eigenvalues that stay separated from the bulk of this matrix's eigenvalues, then $M$ has in probability the same property, and the angle between the corresponding eigenspaces for $K$ and $M$ go to 0 in probability.

(Note that the statements we just made assume that both $M$ and $K$ are symmetric, which is the case here.)

The strategy for the proof is the following. According to the results of Lemma A.3, the matrix $X_i'X_j/p$ has "small" entries off the diagonal, whereas on the diagonal, the entries are essentially constant and equal to $\text{trace}(\Sigma_p)/p$. Hence, it is natural to try to use the $\delta$-method (i.e., do a Taylor expansion) entry by entry. By contrast to standard problems in Statistics, the fact that we have to perform $n^2$ of those Taylor expansions means that the second order term is not negligible a priori. The proof shows that this approach can be carried out rigorously, and that, perhaps surprisingly, the second order term is not too complicated to approximate in operator norm. It is also shown that the third order term plays essentially no role.

Before we start the proof, we want to mention that we will drop the index $p$ in $\Sigma_p$ below to avoid cumbersome notation. Let us also note, more technically, that an important step of the proof is to show that, when the $Y_i$'s have enough moments, they can be treated without much error in spectral results has bounded random variables—the bound depending on the number of moments, $n$ and $p$. This then enables us to use concentration results for convex Lipschitz functions of independent bounded random variables at various important points of the proof and also in Lemma A.3 whose results underly much of the approach taken here.

PROOF OF THEOREM 2.1. First, let us call
$$\tau \triangleq \frac{\text{trace}(\Sigma)}{p}.$$

Using Taylor expansions, we can rewrite our kernel matrix as
$$f(X_i'X_j/p) = f(0) + f'(0)X_i'X_j/p + \frac{f''(0)}{2}(X_i'X_j/p)^2$$
$$+ \frac{f^{(3)}(\xi_{i,j})}{6}(X_i'X_j/p)^3 \quad \text{if } i \neq j,$$
$$f(\|X_i\|_2^2/p) = f(\tau) + f'(\xi_{i,i})\left(\frac{\|X_i\|_2^2}{p} - \tau\right) \quad \text{on the diagonal.}$$

The proof can be separated in different steps. We will break the kernel matrix into a diagonal term and an off-diagonal term. The results of Lemma A.3, after they are shown, will allow us to take care of the diagonal matrix at relatively lost cost. So we postpone that part of the analysis to the end of the proof and we first focus on the off-diagonal matrix.

In what follows, we call "second order term" the matrix $A$ with entries,
$$A_{i,j} = \frac{f''(0)}{2}(X_i'X_j/p)^2 1_{i \neq j}.$$

We call "third order term" the matrix $B$ with entries,

$$B_{i,j} = \frac{f^{(3)}(\xi_{i,j})}{6}(X_i'X_j/p)^3 1_{i \neq j}.$$

The "off-diagonal" matrix is the sum $A + B$.

(A) *Study of the off-diagonal matrix.*

• *Truncation and centralization.* Following the arguments of Lemma 2.2 in [45], we see that because we have assumed that we have $4 + \varepsilon$ absolute moments, and $n \asymp p$, the array $Y = Y_{1 \leq i \leq n, 1 \leq j \leq p}$ is almost surely equal to the array $\widetilde{Y}$ of same dimensions with

$$\widetilde{Y}_{i,j} = Y_{i,j} 1_{|Y_{i,j}| \leq B_p} \qquad \text{where } B_p = p^{1/2-\delta} \text{ and } \delta > 0.$$

We will therefore carry out the analysis on this $\widetilde{Y}$ array. Note that most of the results we will rely on require vectors of i.i.d. entries with mean 0. Of course, $\widetilde{Y}_{i,j}$ has in general a mean different from 0. In other words, if we call $\mu = \mathbf{E}(\widetilde{Y}_{i,j})$, we need to show that we do not lose anything in operator norm by replacing $\widetilde{Y}_i$'s by $U_i$'s with $U_i = \widetilde{Y}_i - \mu 1$. Note that, as seen in Lemma A.3, by plugging in $t = 1/2 - \delta$ in the notation of this lemma, which corresponds to the $4 + \varepsilon$ moment assumption here, we have

$$|\mu| \leq p^{-3/2-\delta}.$$

Now let us call $S$ the matrix $XX'/p$, except that its diagonal is replaced by zeros. From [45], and the fact that $n/p$ stays bounded, we know that $\|XX'/p\|_2 \leq \sigma_1(\Sigma)\|YY'\|_2/p$ stays bounded. Using Lemma A.3, we see that the diagonal of $XX'/p$ stays bounded a.s. in operator norm. Therefore, $\|S\|_2$ is bounded a.s.

Now, as in the proof of Lemma A.3, we have

$$S_{i,j} = \frac{U_i'\Sigma U_j}{p} + \mu\left(\frac{1'\Sigma U_j}{p} + \frac{1'\Sigma U_i}{p}\right) + \mu^2 \frac{1'\Sigma 1}{p} \triangleq \frac{U_i'\Sigma U_j}{p} + R_{i,j} \qquad \text{a.s.}$$

Note that this equality is true a.s. only because it involves replacing $Y$ by $\widetilde{Y}$. The proof of Lemma A.3 shows that

$$|R_{i,j}| \leq \mu 2\sigma_1^{1/2}(\Sigma)(\sigma_1^{1/2}(\Sigma) + p^{-\delta/2}) + \mu^2 \sigma_1(\Sigma) \qquad \text{a.s.}$$

We conclude that, for some constant $C$,

$$\|R\|_F^2 \leq Cn^2\mu^2 \leq Cn^2 p^{-3-2\delta} \to 0 \qquad \text{a.s.}$$

Therefore $\|R\|_2 \to 0$ a.s. In other words, if we call $S_U$ the matrix with $i,j$ entry $U_i'\Sigma U_j/p$ off the diagonal and 0 on the diagonal,

$$\|S - S_U\|_2 \to 0 \qquad \text{a.s.}$$

Now it is a standard result on Hadamard products (see for instance, [10], Problem I.6.13, or [25], Theorems 5.5.1 and 5.5.15) that for two matrices $A$ and $B$, $|||A \circ B|||_2 \leq |||A|||_2 |||B|||_2$. Since the Hadamard product is commutative, we have

$$S \circ S - S_U \circ S_U = (S + S_U) \circ (S - S_U).$$

We conclude that

$$|||S \circ S - S_U \circ S_U|||_2 \leq |||S - S_U|||_2 (|||S|||_2 + |||S_U|||_2) \to 0 \quad \text{a.s.},$$

since $|||S - S_U|||_2 \to 0$ a.s., and $|||S|||_2$ and hence $|||S_U|||_2$ stay bounded, a.s.

The conclusion of this study is that to approximate the second order term in operator norm, it is enough to work with $S_U$ and not $S$, and hence, very importantly, with bounded random variables with zero mean. Further, the proof of Lemma A.3 makes clear that $\sigma_U^2$, the variance of the $U_{i,j}$'s goes to 1, the variance of the $Y_{i,j}$'s, very fast. So if we can approximate the matrix with $(i,j)$-entry $U_i' \Sigma U_j / (p\sigma_U^2)$ consistently in operator norm by a matrix whose operator norm is bounded, this same matrix will constitute an operator norm approximation of $U_i' \Sigma U_j / p$.

In other words, we can assume that, when working with matrices of dimension $n \times p$, the random variables we will be working with have variance 1 without loss of generality and that they have mean 0 and are bounded by $B_p$, $B_p$ depending on $p$ and going to infinity.

• *Control of the second order term.* We now focus on approximating in operator norm the matrix with $(i,j)$th entry,

$$\frac{f''(0)}{2}(X_i'X_j/p)^2 1_{i \neq j}.$$

As we just explained, we assume from now on in all the work concerning the second order term that the vectors $Y_i$ have mean 0, and that their entries have variance 1 and are bounded by $B_p = p^{1/2-\delta}$. This is because we just saw that replacing $Y_i$ by $U_i/\sigma_U$ would not change (a.s. and asymptotically) the operator norm of the matrix to be studied. We note that to make clear that the truncation depends on $p$, we might have wanted to use the notation $Y_i^{(p)}$, but since there will be no ambiguity in the proof, we chose to use the less cumbersome notation $Y_i$.

The control of the second order term turns out to be the most delicate part of the analysis, and the only place where we need the assumption that $X_i = \Sigma^{1/2} Y_i$. Let us call $W$ the matrix with entries

$$W_{i,j} = \begin{cases} \dfrac{(X_i'X_j)^2}{p^2}, & \text{if } i \neq j, \\ 0, & \text{if } i = j. \end{cases}$$

Note that when $i \neq j$,

$$\mathbf{E}(W_{i,j}) = \mathbf{E}(\mathrm{trace}(X_i'X_jX_j'X_i))/p^2 = \mathbf{E}(\mathrm{trace}(X_jX_j'X_iX_i'))/p^2$$
$$= \mathrm{trace}(\Sigma^2)/p^2.$$

Because we assume that $\mathrm{trace}(\Sigma)/p$ has a finite limit, and $n/p$ stays bounded away from 0, we see that the matrix $\mathbf{E}(W)$ has a largest eigenvalue that, in general, does not go to 0. Note also that under our assumptions, $\mathbf{E}(W_{i,j}) = O(1/p)$. Our aim is to show that $W$ can be approximated in operator norm by this constant matrix. So let us consider the matrix $\widetilde{W}$ with entries

$$\widetilde{W}_{i,j} = \begin{cases} \dfrac{(X_i'X_j)^2}{p^2} - \mathrm{trace}(\Sigma^2)/p^2, & \text{if } i \neq j, \\ 0, & \text{if } i = j. \end{cases}$$

Simple computations show that the expected Frobenius norm squared of this matrix does not go to 0. Hence more subtle arguments are needed to control its operator norm. We will show that $\mathbf{E}(\mathrm{trace}(\widetilde{W}^4))$ goes to zero which implies that $\mathbf{E}(\|\widetilde{W}\|_2^4)$ goes to zero because $\widetilde{W}$ is real symmetric.

The elements contributing to $\mathrm{trace}(\widetilde{W}^4)$ are generally of the form $\widetilde{W}_{i,j}\widetilde{W}_{j,k} \times \widetilde{W}_{k,l}\widetilde{W}_{l,i}$. We are going to study these terms according to how many indices are equal to each other.

(i) *Terms involving 4 different indices: $i \neq j \neq k \neq l$.* We first focus on the case where all these indices $(i,j,k,l)$ are different. Recall that $X_i = \Sigma^{1/2}Y_i$, where $Y_i$ has i.i.d. entries. We want to compute $\mathbf{E}(\widetilde{W}_{i,j}\widetilde{W}_{j,k}\widetilde{W}_{k,l}\widetilde{W}_{l,i})$, so it is natural to focus first on

$$\mathbf{E}(\widetilde{W}_{i,j}\widetilde{W}_{j,k}\widetilde{W}_{k,l}\widetilde{W}_{l,i}|Y_i,Y_k).$$

Now, note that

$$\widetilde{W}_{i,j} = \frac{1}{p^2}\{Y_i'\Sigma Y_j Y_j'\Sigma Y_i - \mathrm{trace}(\Sigma^2)\}$$
$$= \frac{1}{p^2}\{Y_i'\Sigma(Y_jY_j' - \mathrm{Id})\Sigma Y_i + \mathrm{trace}(\Sigma^2(Y_iY_i' - \mathrm{Id}))\}.$$

Hence, calling

$$M_j \triangleq Y_jY_j' - \mathrm{Id},$$

we have

$$p^4\widetilde{W}_{i,j}\widetilde{W}_{j,k} = (Y_i'\Sigma M_j \Sigma Y_i Y_k'\Sigma M_j \Sigma Y_k) + (Y_i'\Sigma M_j \Sigma Y_i)\mathrm{trace}(\Sigma^2 M_k)$$
$$+ (Y_k'\Sigma M_j \Sigma Y_k)\mathrm{trace}(\Sigma^2 M_i) + \mathrm{trace}(\Sigma^2 M_i)\mathrm{trace}(\Sigma^2 M_k).$$

Now, of course, we have $\mathbf{E}(M_j) = \mathbf{E}(M_j|Y_i, Y_k) = 0$. Hence,

$$p^4 \mathbf{E}(\widetilde{W}_{i,j}\widetilde{W}_{j,k}|Y_i, Y_k) = (Y_i' \Sigma \mathbf{E}(M_j \Sigma Y_i Y_k' \Sigma M_j|Y_i, Y_k) \Sigma Y_k)$$
$$+ \operatorname{trace}(\Sigma^2 M_i) \operatorname{trace}(\Sigma^2 M_k).$$

If $M$ is a deterministic matrix, we have, since $\mathbf{E}(Y_j Y_j') = \mathrm{Id}$,

$$\mathbf{E}(M_j M M_j) = \mathbf{E}(Y_j Y_j' M Y_j Y_j') - M.$$

If we now use Lemma A.1, and, in particular, (4), page 41, we finally have, recalling that here $\sigma^2 = 1$,

$$\mathbf{E}(M_j M M_j) = (M + M') + (\mu_4 - 3) \operatorname{diag}(M) + \operatorname{trace}(M) \operatorname{Id} - M$$
$$= M' + (\mu_4 - 3) \operatorname{diag}(M) + \operatorname{trace}(M) \operatorname{Id}.$$

In the case of interest here, we have $M = \Sigma Y_i Y_k' \Sigma$, and the expectation is to be understood conditionally on $Y_i, Y_k$, but because we have assumed that the indices are different and the $Y_m$'s are independent, we can do the computation of the conditional expectation as if $M$ were deterministic. Therefore, we have

$$(Y_i' \Sigma \mathbf{E}(M_j \Sigma Y_i Y_k' \Sigma M_j|Y_i, Y_k) \Sigma Y_k)$$
$$= Y_i' \Sigma [\Sigma Y_k Y_i' \Sigma + (\mu_4 - 3) \operatorname{diag}(\Sigma Y_i Y_k' \Sigma) + (Y_k' \Sigma^2 Y_i) \operatorname{Id}] \Sigma Y_k$$
$$= [(Y_i' \Sigma^2 Y_k)^2 + (\mu_4 - 3) Y_i' \Sigma \operatorname{diag}(\Sigma Y_i Y_k' \Sigma) \Sigma Y_k + (Y_i' \Sigma^2 Y_k)^2].$$

Naturally, we have $\mathbf{E}(\widetilde{W}_{i,j}\widetilde{W}_{j,k}|Y_i, Y_k) = \mathbf{E}(\widetilde{W}_{k,l}\widetilde{W}_{l,i}|Y_i, Y_k)$, and therefore, by using properties of conditional expectation, since all the indices are different,

$$p^8 \mathbf{E}(\widetilde{W}_{i,j}\widetilde{W}_{j,k}\widetilde{W}_{k,l}\widetilde{W}_{l,i})$$
$$= \mathbf{E}([2(Y_i' \Sigma^2 Y_k)^2 + (\mu_4 - 3) Y_i' \Sigma \operatorname{diag}(\Sigma Y_i Y_k' \Sigma) \Sigma Y_k$$
$$+ \operatorname{trace}(\Sigma^2 M_i) \operatorname{trace}(\Sigma^2 M_k)]^2).$$

By convexity, we have $(a+b+c)^2 \leq 3(a^2+b^2+c^2)$, so to control the above expression, we just need to control the square of each of the terms appearing in it. In other words, we need to understand the terms

$$T_1 = \mathbf{E}((Y_i' \Sigma^2 Y_k)^4),$$
$$T_2 = \mathbf{E}([Y_i' \Sigma \operatorname{diag}(\Sigma Y_i Y_k' \Sigma) \Sigma Y_k]^2)$$

and

$$T_3 = \mathbf{E}([\operatorname{trace}(\Sigma^2 M_i) \operatorname{trace}(\Sigma^2 M_k)]^2).$$

*Study of $T_1$.* Let us start by the term $T_1 = \mathbf{E}((Y_i' \Sigma^2 Y_k)^4)$. A simple rewriting shows that

$$(Y_i' \Sigma^2 Y_k)^4 = Y_i' \Sigma^2 Y_k Y_k' \Sigma^2 Y_i Y_i' \Sigma^2 Y_k Y_k' \Sigma^2 Y_i.$$

Using (4) in Lemma A.1, we therefore have, using the fact that $\Sigma^2 Y_i Y_i' \Sigma^2$ is symmetric,

$$\begin{aligned}
\mathbf{E}((Y_i' \Sigma^2 Y_k)^4 | Y_i) \\
&= Y_i' \Sigma^2 [2\Sigma^2 Y_i Y_i' \Sigma^2 + (\mu_4 - 3) \operatorname{diag}(\Sigma^2 Y_i Y_i' \Sigma^2) \\
&\qquad\qquad + \operatorname{trace}(\Sigma^2 Y_i Y_i' \Sigma^2) \operatorname{Id}] \Sigma^2 Y_i \\
&= 3(Y_i' \Sigma^4 Y_i)^2 + (\mu_4 - 3) Y_i' \Sigma^2 \operatorname{diag}(\Sigma^2 Y_i Y_i' \Sigma^2) \Sigma^2 Y_i.
\end{aligned}$$

Finally, we have, using (5) in Lemma A.1,

$$\begin{aligned}
\mathbf{E}((Y_i' \Sigma^2 Y_k)^4) &= 3[2 \operatorname{trace}(\Sigma^4) + (\operatorname{trace}(\Sigma^4))^2 + (\mu_4 - 3) \operatorname{trace}(\Sigma^4 \circ \Sigma^4)] \\
&\quad + (\mu_4 - 3) \mathbf{E}(Y_i' \Sigma^2 \operatorname{diag}(\Sigma^2 Y_i Y_i' \Sigma^2) \Sigma^2 Y_i).
\end{aligned}$$

Now we have

$$\begin{aligned}
Y_i' \Sigma^2 \operatorname{diag}(\Sigma^2 Y_i Y_i' \Sigma^2) \Sigma^2 Y_i &= \operatorname{trace}(\Sigma^2 Y_i Y_i' \Sigma^2 \operatorname{diag}(\Sigma^2 Y_i Y_i' \Sigma^2)) \\
&= \operatorname{trace}(\Sigma^2 Y_i Y_i' \Sigma^2 \circ \Sigma^2 Y_i Y_i' \Sigma^2).
\end{aligned}$$

Calling $v_i = \Sigma^2 Y_i$, we note that the matrix whose trace is taken is $(v_i v_i') \circ (v_i v_i') = (v_i \circ v_i)(v_i \circ v_i)'$ (see [24], page 458 or [25], page 307). Hence,

$$Y_i' \Sigma^2 \operatorname{diag}(\Sigma^2 Y_i Y_i' \Sigma^2) \Sigma^2 Y_i = \|v_i \circ v_i\|_2^2.$$

Now let us call $m_k$ the $k$th column of the matrix $\Sigma^2$. Using the fact that $\Sigma^2$ is symmetric, we see that the $k$th entry of the vector $v_i$ is $v_i(k) = m_k' Y_i$. So $v_i(k)^4 = Y_i' m_k m_k' Y_i Y_i' m_k m_k' Y_i$. Calling $\mathcal{M}_k = m_k m_k'$, we see, using (5) in Lemma A.1, that

$$\mathbf{E}(v_i(k)^4) = 2 \operatorname{trace}(\mathcal{M}_k^2) + [\operatorname{trace}(\mathcal{M}_k)]^2 + (\mu_4 - 3) \operatorname{trace}(\mathcal{M}_k \circ \mathcal{M}_k).$$

Using the definition of $\mathcal{M}_k$, we finally get that

$$\mathbf{E}(v_i(k)^4) = 3\|m_k\|_2^4 + (\mu_4 - 3)\|m_k \circ m_k\|_2^2.$$

Now, note that if $C$ is a generic matrix and $C_k$ is its $k$th column, denoting by $e_k$ the $k$th vector of the canonical basis, we have $C_k = Ce_k$, and hence $\|C_k\|_2^2 = e_k' C' C e_k \leq \sigma_1^2(C)$ where $\sigma_1(C)$ is the largest singular value of $C$. So in particular, if we call $\lambda_1(D)$ the largest eigenvalue of a positive semi-definite matrix $D$, we have $\|m_k\|_2^4 \leq \lambda_1(\Sigma^4)\|m_k\|_2^2$.

After recalling the definition of $m_k$, and using the fact that $\sum_k \|m_k \circ m_k\|_2^2 = \|\Sigma^2 \circ \Sigma^2\|_F^2$, we deduce that

$$\begin{aligned}
\mathbf{E}(\|v_i \circ v_i\|_2^2) &= 3 \sum_k \|m_k\|_2^4 + (\mu_4 - 3) \sum_k \|m_k \circ m_k\|_2^2 \\
&\leq 3\lambda_1(\Sigma^4) \operatorname{trace}(\Sigma^4) + (\mu_4 - 3) \operatorname{trace}([\Sigma^2 \circ \Sigma^2]^2).
\end{aligned}$$

Therefore, we can conclude that

$$\mathbf{E}((Y_i'\Sigma^2 Y_k)^4) \leq 3\lambda_1(\Sigma^4)\operatorname{trace}(\Sigma^4) + (\mu_4 - 3)\operatorname{trace}([\Sigma^2 \circ \Sigma^2]^2).$$

Now recall that, according to Theorem 5.5.19 in [25], if $C$ and $D$ are positive semidefinite matrices, $\lambda(C \circ D) \prec_w d(C) \circ \lambda(D)$ where $\lambda(D)$ is the vector of decreasingly ordered eigenvalues of $D$, and $d(C)$ denotes the vector of decreasingly ordered diagonal entries of $C$ (because all the matrices are positive semidefinite, their eigenvalues are their singular values). Here $\prec_w$ denotes weak (sub)majorization. In our case, of course, $C = D = \Sigma^2$. Using the results of Example II.3.5(iii) in [10], with the function $\phi(x) = x^2$, we see that

$$\operatorname{trace}((\Sigma^2 \circ \Sigma^2)^2) = \sum \lambda_i^2(\Sigma^2 \circ \Sigma^2) \leq \sum d_i^2(\Sigma^2)\lambda_i^2(\Sigma^2) \leq \lambda_1(\Sigma^4)\operatorname{trace}(\Sigma^4).$$

Finally, we have

$$T_1 = \mathbf{E}((Y_i'\Sigma^2 Y_k)^4) \leq (3 + |\mu_4 - 3|)\lambda_1(\Sigma^4)\operatorname{trace}(\Sigma^4). \tag{1}$$

This bounds the first term, $T_1$, in our upper bound.

*Study of $T_3$.* Let us now turn to the third term, $T_3 = \mathbf{E}([\operatorname{trace}(\Sigma^2 M_i)\operatorname{trace} \times (\Sigma^2 M_k)]^2)$. We remind the reader that $M_i = Y_i Y_i' - \operatorname{Id}$. By independence of $Y_i$ and $Y_k$, it is enough to understand $\mathbf{E}([\operatorname{trace}(\Sigma^2 M_i)]^2)$. Note that

$$\mathbf{E}([\operatorname{trace}(\Sigma^2 M_i)]^2) = \mathbf{E}([Y_i'\Sigma^2 Y_i - \operatorname{trace}(\Sigma^2)]^2)$$
$$= \mathbf{E}(Y_i'\Sigma^2 Y_i Y_i'\Sigma^2 Y_i) - \operatorname{trace}(\Sigma^2)^2.$$

Using (5) in Lemma A.1, we conclude that

$$\mathbf{E}([\operatorname{trace}(\Sigma^2 M_i)]^2) = 2\operatorname{trace}(\Sigma^4) + (\mu_4 - 3)\operatorname{trace}(\Sigma^2 \circ \Sigma^2).$$

Using the fact that we know the diagonal of $\Sigma^2 \circ \Sigma^2$, we conclude that

$$T_3 = \mathbf{E}([\operatorname{trace}(\Sigma^2 M_i)]^2 [\operatorname{trace}(\Sigma^2 M_k)]^2)$$
$$\leq \{2\operatorname{trace}(\Sigma^4) + |\mu_4 - 3|\lambda_1(\Sigma^2)\operatorname{trace}(\Sigma^2)\}^2. \tag{2}$$

So we have an upper bound on $T_3$.

*Study of $T_2$.* Finally, let us turn to the middle term, $T_2 = \mathbf{E}([Y_i'\Sigma\operatorname{diag}(\Sigma Y_i Y_k' \times \Sigma)\Sigma Y_k]^2)$. Before we square it, the argument of the expectation has the form $Y_i'\Sigma\operatorname{diag}(\Sigma Y_k Y_i'\Sigma)\Sigma Y_k$. Call $u_k = \Sigma Y_k$. Making the same computations as above, we find that

$$Y_i'\Sigma\operatorname{diag}(\Sigma Y_k Y_i'\Sigma)\Sigma Y_k = \operatorname{trace}(\operatorname{diag}(\Sigma Y_k Y_i'\Sigma)\Sigma Y_k Y_i'\Sigma)$$
$$= \operatorname{trace}((\Sigma Y_k Y_i'\Sigma) \circ (\Sigma Y_k Y_i'\Sigma))$$
$$= \operatorname{trace}((u_k u_i') \circ (u_k u_i')) = \operatorname{trace}((u_k \circ u_k)(u_i \circ u_i)')$$
$$= (u_i \circ u_i)'(u_k \circ u_k).$$

We deduce, using independence and elementary properties of inner products, that

$$\mathbf{E}([Y_i'\Sigma\operatorname{diag}(\Sigma Y_k Y_i'\Sigma)\Sigma Y_k]^2) \leq \mathbf{E}(\|u_i \circ u_i\|_2^2)\mathbf{E}(\|u_k \circ u_k\|_2^2).$$

Note that to arrive at (1), we studied expressions similar to $\mathbf{E}(\|u_i \circ u_i\|_2^2)$. So we can similarly conclude that

$$T_2 = \mathbf{E}([Y_i'\Sigma\operatorname{diag}(\Sigma Y_k Y_i'\Sigma)\Sigma Y_k]^2) \leq \{(3+|\mu_4-3|)\lambda_1(\Sigma^2)\operatorname{trace}(\Sigma^2)\}^2. \tag{3}$$

With our assumptions, the terms (1), (2) and (3) are $\mathrm{O}(p^2)$. Note that in the computation of the trace, there are $\mathrm{O}(n^4)$ such terms. Finally, note that the expectation of interest to us corresponds to the sum of the three quadratic terms divided by $p^8$. So the total contribution of these terms is in expectation $\mathrm{O}(p^{-2})$. This takes care of the contribution of the terms involving four different indices, as it shows that

$$0 \leq \mathbf{E}\bigg(\sum_{i\neq j\neq k\neq l} \widetilde{W}_{i,j}\widetilde{W}_{j,k}\widetilde{W}_{k,l}\widetilde{W}_{l,i}\bigg) = \mathrm{O}(p^{-2}).$$

(ii) *Terms involving three different indices: $i \neq j \neq k$.* Note that because $\widetilde{W}_{i,i} = 0$, terms involving 3 different indices with a nonzero contribution are necessarily of the form $(\widetilde{W}_{i,j})^2(\widetilde{W}_{i,k})^2$, since terms with a cycle of length 3 all involve a term of the form $\widetilde{W}_{i,i}$ and hence contribute 0. Let us now focus on those terms, assuming that $j \neq k$. Note that we have $\mathrm{O}(n^3)$ such terms and that it is enough to focus on the $W_{i,j}^2 W_{i,k}^2$, since the contribution of the other terms is, in expectation, of order $1/p^4$ [with our assumptions $\operatorname{trace}(\Sigma^2)/p^2 = \mathrm{O}(1/p)$], and because we have only $n^3$ terms in the sum, this extra contribution is asymptotically zero. Now, we clearly have $\mathbf{E}(W_{i,j}^2 W_{i,k}^2|Y_i) = [\mathbf{E}(W_{i,j}^2|Y_i)]^2$, by conditional independence of the two terms. The computation of $\mathbf{E}(W_{i,j}^2|Y_i)$ is similar to the ones we have made above, and we have

$$p^4\mathbf{E}(W_{i,j}^2|Y_i) = 2(Y_i'\Sigma^2 Y_i)^2 + (\mu_4-3)Y_i'\Sigma\operatorname{diag}(\Sigma Y_i Y_i'\Sigma)\Sigma Y_i$$
$$+ (\operatorname{trace}(\Sigma Y_i Y_i'\Sigma))^2.$$

Using the fact that $\mathcal{K}_i = \Sigma Y_i Y_i'\Sigma$ is positive semidefinite, and hence its diagonal entries are nonnegative, we have $\operatorname{trace}(\mathcal{K}_i \circ \mathcal{K}_i) \leq (\operatorname{trace}(\mathcal{K}_i))^2$, and we conclude that

$$p^4\mathbf{E}(W_{i,j}^2|Y_i) \leq (3+|\kappa_4-3|)(Y_i'\Sigma^2 Y_i)^2 \leq (3+|\kappa_4-3|)\sigma_1(\Sigma)^4\|Y_i\|_2^4.$$

Hence,

$$\mathbf{E}(W_{i,j}^2 W_{i,k}^2) \leq \frac{1}{p^8}(3+|\kappa_4-3|)^2\sigma_1(\Sigma)^8\|Y_i\|_2^8.$$

Now, the application $F$ which takes a vector and returns its Euclidean norm is trivially a convex 1-Lipschitz function, with respect to Euclidean norm. Because the entries of $Y_i$ are bounded by $B_p$, we see that, according to Corollary 4.10 in [29], $F(Y_i) = \|Y_i\|_2$ satisfies a concentration inequality, namely, for $r > 0$, $P(|\|Y_i\|_2 - m_F| > r) \leq 4\exp(-r^2/16B_p^2)$ where $m_F$ is a median of $F(Y_i) = \|Y_i\|_2$ (hence $m_F$ is a deterministic quantity). A simple integration (see, for instance, the proof of Proposition 1.9 in [29], and change the power from 2 to 8) then shows that

$$\mathbf{E}(|\|Y_i\|_2 - m_F|^8) = \mathrm{O}(B_p^8).$$

Now we know, according to Proposition 1.9 in [29], that if $\mu_F$ is the mean of $F(Y_i)$, that is, $\mu_F = \mathbf{E}(\|Y_i\|_2)$, $\mu_F$ exists and $|m_F - \mu_F| = \mathrm{O}(B_p)$. Since $\mu_F^2 \leq \mu_{F^2} = \mathbf{E}(\|Y_i\|_2^2) = p$, we conclude that, if $C$ denotes a generic constant that may change from display to display,

$$\mathbf{E}(\|Y_i\|_2^8) \leq \mathbf{E}(|\|Y_i\|_2 - m_F + m_F|^8) \leq 2^7(\mathbf{E}(|\|Y_i\|_2 - m_F|^8) + m_F^8)$$
$$\leq C(\mathbf{E}(|\|Y_i\|_2 - m_F|^8) + |m_F - \mu_F|^8 + \mu_F^8) \leq C(B_p^8 + p^4).$$

Now our original assumption about the number of absolute moments of the random variables of interest imply that $B_p = \mathrm{O}(p^{1/2-\delta})$. Consequently,

$$\mathbf{E}(\|Y_i\|_2^8) = \mathrm{O}(p^4).$$

Therefore,

$$\mathbf{E}(W_{i,j}^2 W_{i,k}^2) = \mathrm{O}(p^{-4})$$

and

$$\sum_i \sum_{j \neq i, k \neq i, j \neq k} \mathbf{E}(W_{i,j}^2 W_{i,k}^2) = \mathrm{O}(p^{-1}).$$

Hence, we also have

$$\sum_i \sum_{j \neq i, k \neq i, j \neq k} \mathbf{E}(\widetilde{W}_{i,j}^2 \widetilde{W}_{i,k}^2) = \mathrm{O}(p^{-1}).$$

(iii) *Terms involving two different indices: $i \neq j$.* The last terms we have to focus on to control $\mathbf{E}(\mathrm{trace}(\widetilde{W}^4))$ are of the form $\widetilde{W}_{i,j}^4$. Note that we have $n^2$ terms like this. Since by convexity, $(a+b)^4 \leq 8(a^4 + b^4)$, we see that it is enough to understand the contribution of $W_{i,j}^4$ to show that $\sum_{i,j} \mathbf{E}(\widetilde{W}_{i,j}^4)$ tends to zero. Now, let us call for a moment $v = \Sigma Y_i$ and $u = Y_j$. The quantity of interest to us is basically of the form $\mathbf{E}((u'v)^8)$. Let us do computations conditional on $v$. We note that since the entries of $u$ are independent and have mean 0, in the expansion of $(u'v)^8$, the only terms that will contribute a nonzero quantity to the expectation have entries of $u$ raised to a power greater than 2. We can decompose the sum representing $\mathbf{E}((u'v)^8|v)$ into subterms, according to what powers of the terms are involved. There are 6 terms: $(2,2,2,2)$ (i.e., all terms are raised to the power 2), $(3,3,2)$ (i.e., two terms are raised to the power 3, and one to the power 2), $(4,2,2)$, $(4,4)$, $(5,3)$, $(6,2)$ and $(8)$. For instance the subterm corresponding to $(2,2,2,2)$ is, before taking expectations,

$$\sum_{i_1 \neq i_2 \neq i_3 \neq i_4} u_{i_1}^2 u_{i_2}^2 u_{i_3}^2 u_{i_4}^2 (v_{i_1} v_{i_2} v_{i_3} v_{i_4})^2.$$

After taking expectations conditional on $v$, we see that it is obviously non-negative and contributes

$$(\sigma^2)^4 \sum_{i_1 \neq i_2 \neq i_3 \neq i_4} (v_{i_1} v_{i_2} v_{i_3} v_{i_4})^2 \leq \left(\sum v_i^2\right)^4 = (Y_i' \Sigma^2 Y_i)^4 \leq \sigma_1(\Sigma)^8 \|Y_i\|_2^8.$$

Note that we just saw that $\mathbf{E}(\|Y_i\|_2^8) = \mathrm{O}(p^4)$ in our context. Similarly, the term $(3,3,2)$ will contribute

$$\mu_3^2 \sigma^2 \sum_{i_1 \neq i_2 \neq i_3} v_{i_1}^3 v_{i_2}^3 v_{i_3}^2.$$

In absolute value, this term is less than

$$\mu_3^2 \sigma^2 \left(\sum |v_i|^3\right)^2 \left(\sum v_i^2\right).$$

Now, note that if $z$ is such that $\|z\|_2 = 1$, we have, for $p \geq 2$, $\sum |z_i|^p \leq \sum z_i^2 = 1$. Applied to $z = v/\|v\|_2$, we conclude that $\sum |v_i|^p \leq \|v\|_2^p$. Consequently, the term $(3,3,2)$ contributes in absolute value less than

$$\mu_3^2 \sigma^2 \|v\|_2^8.$$

The same analysis can be repeated for all the other terms which are all found to be less than $\|v\|_2^8$ times the moments of $u$ involved. Because we have assumed that our original random variables had $4 + \varepsilon$ absolute moments, the moments of order less than 4 cause no problem. The moments of order higher than 4, say $4 + k$, can be bounded by $\mu_4 B_p^k$. Consequently, we see that

$$\mathbf{E}(W_{i,j}^4) = \mathbf{E}(\mathbf{E}(W_{i,j}^4|Y_i)) \leq C B_p^4 \mathbf{E}\left(\frac{\|Y_i\|^8}{p^8}\right) = \mathrm{O}(B_p^4/p^4) = \mathrm{O}(p^{-(2+4\delta)}).$$

Since we have $n^2$ such terms, we see that

$$\sum_{i \neq j} \mathbf{E}(W_{i,j}^4) \to 0 \quad \text{as } p \to \infty.$$

Using our earlier convexity remark, we finally conclude that

$$\sum_{i \neq j} \mathbf{E}(\widetilde{W}_{i,j}^4) \to 0 \quad \text{as } p \to \infty.$$

(iv) *Second order term: combining all the elements.* We have therefore established control of the second order term and seen that the largest singular value of $\widetilde{W}$ goes to 0 in probability, using Chebyshev's inequality. Note that we have also shown that the operator norm of $W$ is bounded in probability and that

$$\left\|\!\left\| W - \frac{\operatorname{trace}(\Sigma^2)}{p^2}(11' - \operatorname{Id}) \right\|\!\right\|_2 \to 0 \quad \text{in probability.}$$

• *Control of the third order term.* We note that the third order term is of the form $f^{(3)}(\xi_{i,j})\frac{X_i'X_j}{p}W_{i,j}$. According to Lemma A.5, if $M$ is a real symmetric matrix with nonnegative entries, and $E$ is a symmetric matrix such that $\max_{i,j}|E_{i,j}| = \zeta$, then

$$\sigma_1(E \circ M) \leq \zeta \sigma_1(M).$$

Note that $W$ is real symmetric matrix with nonnegative entries. So all we have to show to prove that the third order term goes to zero in operator norm is that $\max_{i \neq j}|X_i'X_j/p|$ goes to 0 because we have just established that $\|\!\|W\|\!\|_2$ remains bounded in probability. We are going to make use of Lemma A.3, page 45 in the Appendix. In our setting, we have $B_p = p^{1/2-\delta}$, or $2/m = 1/2 - \delta$. The lemma implies, for instance, that

$$\max_{i \neq j}|X_i'X_j/p| \leq p^{-\delta}\log(p) \quad \text{a.s.}$$

So $\max_{i \neq j}|X_i'X_j/p| \to 0$ a.s. Note that this implies that $\max_{i \neq j}|\xi_{i,j}| \to 0$ a.s. Since we have assumed that $f^{(3)}$ exists and is continuous and hence bounded in a neighborhood of 0, we conclude that

$$\max_{i,j}|f^{(3)}(\xi_{i,j})X_i'X_j/p| = \operatorname{o}(p^{-\delta/2}) \quad \text{a.s.}$$

If we call $E$ the matrix with entry $E_{i,j} = f^{(3)}(\xi_{i,j})X_i'X_j/p$ off-the diagonal and 0 on the diagonal, we see that $E$ satisfies the conditions put forth in our discussion earlier in this section and we conclude that

$$\|\!\|E \circ W\|\!\|_2 \leq \max_{i,j}|E_{i,j}|\|\!\|W\|\!\|_2 = \operatorname{o}(p^{-\delta/2}) \quad \text{a.s.}$$

Hence, the operator norm of the third order term goes to 0 almost surely. [To maybe clarify our arguments, let us repeat that we analyzed the second order term by replacing the $Y_i$'s by, in the notation of the truncation and centralization discussion, $U_i$. Let us call $W_U = S_U \circ S_U$, again using notation introduced in the truncation and centralization discussion. As we saw, $\|\!\|W - W_U\|\!\|_2 \to 0$ a.s., so showing, as we did, that $\|\!\|W_U\|\!\|_2$ remains bounded (a.s.) implies that $\|\!\|W\|\!\|_2$ does too, and this is the only thing we need in our argument showing the control of the third order term.]

(B) *Control of the diagonal term.* The proof here is divided into two parts. First, we show that the error term coming from the first order expansion of the diagonal is easily controlled. Then we show that the terms added when replacing the off-diagonal matrix by $XX'/p + \text{trace}(\Sigma^2)/p^2 11'$ can also be controlled. Recall the notation $\tau = \text{trace}(\Sigma)/p$.

• *Errors induced by diagonal approximation.* Note that Lemma A.3 guarantees that for all $i$, $|\xi_{i,i} - \tau| \leq p^{-\delta/2}$, a.s. Because we have assumed that $f'$ is continuous and hence bounded in a neighborhood of $\tau$, we conclude that $f'(\xi_{i,i})$ is uniformly bounded in $p$. Now Lemma A.3 also guarantees that

$$\max_i \left| \frac{\|X_i\|_2^2}{p} - \tau \right| \leq p^{-\delta} \qquad \text{a.s.}$$

Hence, the diagonal matrix with entries $f(\|X_i\|_2^2/p)$ can be approximated consistently in operator norm by $f(\tau)\,\text{Id}$ a.s.

• *Errors induced by off-diagonal approximation.* When we replace the off-diagonal matrix by $f'(0)XX'/p + [f(0) + f''(0)\text{trace}(\Sigma^2)/2p^2]11'$, we add a diagonal matrix with $(i,i)$ entry $f(0) + f'(0)\|X_i\|_2^2/p + f''(0)\text{trace}(\Sigma^2)/2p^2$ which we need to subtract eventually. We note that $0 \leq \text{trace}(\Sigma^2)/p^2 \leq \sigma_1^2(\Sigma)/p \to 0$ when $\sigma_1(\Sigma)$ remains bounded in $p$. So this term does not create any problem. Now, we just saw that the diagonal matrix with entries $\|X_i\|_2^2/p$ can be consistently approximated in operator norm by $(\text{trace}(\Sigma)/p)\,\text{Id}$. So the diagonal matrix with $(i,i)$ entry $f(0) + f'(0)\|X_i\|_2^2/p + f''(0)\text{trace}(\Sigma^2)/2p^2$ can be approximated consistently in operator norm by $(f(0) + f'(0)\text{trace} \times (\Sigma)/p)\,\text{Id}$ a.s.

This finishes the proof. □

2.3. *Kernel random matrices of the type $f(\|X_i - X_j\|_2^2/p)$.* As is to be expected, the properties of such matrices can be deduced from the study of inner product kernel matrices, with a little bit of extra work. We need to slightly modify the distributional assumptions under which we work, and consider the case where we have $5 + \varepsilon$ absolute moments for the entries of $Y_i$. We also need to assume that $f$ is regular is the neighborhood of different points. Otherwise, the assumptions are the same as that of Theorem 2.1. We have the following theorem:

THEOREM 2.2 (Spectrum of Euclidean distance kernel matrices). *Consider the $n \times n$ kernel matrix $M$ with entries*

$$M_{i,j} = f\left(\frac{\|X_i - X_j\|_2^2}{p}\right).$$

*Let us call*

$$\tau = 2\frac{\text{trace}(\Sigma)}{p}.$$

*Let us call $\psi$ the vector with $i$th entry $\psi_i = \|X_i\|_2^2/p - \text{trace}(\Sigma)/p$. Suppose that the assumptions of Theorem 2.1 hold, but that conditions* (e) *and* (f) *are replaced by:*

(e′) *The entries of $Y_i$, a $p$-dimensional random vector, are i.i.d. Also, denoting by $Y_i(k)$ the $k$th entry of $Y_i$, we assume that $\mathbf{E}(Y_i(k)) = 0$, $\text{var}(Y_i(k)) = 1$ and $\mathbf{E}(|Y_i(k)|^{5+\varepsilon}) < \infty$ for some $\varepsilon > 0$. (We say that $Y_i$ has $5 + \varepsilon$ absolute moments.)*

(f′) *$f$ is $C^3$ in a neighborhood of $\tau$.*

*Then $M$ can be approximated consistently in operator norm (and in probability) by the matrix $K$, defined by*

$$K = f(\tau)11' + f'(\tau)\left[1\psi' + \psi 1' - 2\frac{XX'}{p}\right]$$
$$+ \frac{f''(\tau)}{2}\left[1(\psi \circ \psi)' + (\psi \circ \psi)1' + 2\psi\psi' + 4\frac{\text{trace}(\Sigma^2)}{p^2}11'\right] + v_p\,\text{Id},$$
$$v_p = f(0) + \tau f'(\tau) - f(\tau).$$

*In other words,*

$$|\!|\!|M - K|\!|\!|_2 \to 0 \qquad \text{in probability.}$$

PROOF. Note that here the diagonal is just $f(0)\,\text{Id}$ and it will cause no trouble. The work, therefore, focuses on the off-diagonal matrix. In what follows, we call $\tau = 2\frac{\text{trace}(\Sigma)}{p}$. Let us define

$$A_{i,j} = \frac{\|X_i\|_2^2}{p} + \frac{\|X_j\|_2^2}{p} - \tau$$

and

$$S_{i,j} = \frac{X_i'X_j}{p}.$$

With these notation, we have, off the diagonal, that is, when $i \neq j$, by a Taylor expansion,

$$M_{i,j} = f(\tau) + [A_{i,j} - 2S_{i,j}]f'(\tau) + \frac{1}{2}[A_{i,j} - 2S_{i,j}]^2 f''(\tau)$$
$$+ \frac{1}{6}f^{(3)}(\xi_{i,j})[A_{i,j} - 2S_{i,j}]^3.$$

We note that the matrix $A$ with entries $A_{i,j}$ is a rank 2 matrix. As a matter of fact, it can be written, if $\psi$ is the vector with entries $\psi_i = \frac{\|X_i\|_2^2}{p} - \tau/2$, $A = 1\psi' + \psi 1'$. Using the well-known identity (see, e.g., [23], Chapter 1, Theorem 3.2),

$$\det(I + uv' + vu') = \det\begin{pmatrix} 1 + u'v & \|u\|_2^2 \\ \|v\|_2^2 & 1 + u'v \end{pmatrix};$$

we see immediately that the nonzero eigenvalues of $A$ are

$$1'\psi \pm \sqrt{n}\|\psi\|_2.$$

After these preliminary remarks, we are ready to start the proof per se.

• *Truncation and centralization.* Since we assume $5 + \varepsilon$ absolute moments, we see, using Lemma 2.2 in [45], that we can truncate the $Y_i$'s at level $B_p = p^{2/5-\delta}$ with $\delta > 0$ and a.s. not change the data matrix. We then need to centralize the vectors truncated at $p^{2/5-\delta}$. Note that because we work with $X_i - X_j = \Sigma^{1/2}(Y_i - Y_j)$, centralization creates absolutely no problem here since it is absorbed in the difference. So in what follows we can assume without loss of generality that we are working with vectors $X_i = \Sigma^{1/2}Y_i$ where the entries of $Y_i$ are bounded by $p^{2/5-\delta}$ and $\mathbf{E}(Y_i) = 0$. The issue of variance 1 is addressed as before, so we can assume that the entries of $Y_i$ have variance 1.

• *Concentration of $\|X_i - X_j\|_2^2/p$.* By plugging in the results of Corollary A.2, with $2/m = 2/5 - \delta$, we get that

$$\max_{i \neq j}\left|\frac{\|X_i - X_j\|_2^2}{p} - 2\frac{\text{trace}(\Sigma)}{p}\right| \leq \log(p)p^{-1/10-\delta}.$$

Also, using the result of Lemma A.3, we have

$$\max_i |\psi_i| = \max_i\left|\frac{\|X_i\|_2^2}{p} - \frac{\text{trace}(\Sigma)}{p}\right| \leq \log(p)p^{-1/10-\delta}.$$

Note that, as explained in the proof of Lemma A.3, these results are true whether we work with $Y_i$ or their truncated and centralized version.

• *Control of the second order term.* The second order term is the matrix with $(i,j)$-entry

$$1_{i \neq j}\tfrac{1}{2}f''(\tau)(A_{i,j} - S_{i,j})^2.$$

Let us call $T$ the matrix with 0 on the diagonal and off-diagonal entries $T_{i,j} = (A_{i,j} - 2S_{i,j})^2$. In other words,

$$T_{i,j} = 1_{i \neq j}\left(\frac{\|X_i - X_j\|_2^2 - 2\,\text{trace}(\Sigma)}{p}\right)^2.$$

We simply write $(A_{i,j} - 2S_{i,j})^2 = A_{i,j}^2 - 4A_{i,j}S_{i,j} + 4S_{i,j}^2$. In the notation of the proof of Theorem 2.1, the matrix with entries $S_{i,j}^2$ off the diagonal and 0 on the diagonal is what we called $W$. We have already shown that

$$\left\|\!\left\| W - \frac{\text{trace}(\Sigma^2)}{p^2}(11' - \text{Id}) \right\|\!\right\|_2 \to 0 \quad \text{in probability.}$$

Now, let us focus on the term $A_{i,j}S_{i,j}$. Let us call $H$ the matrix with

$$H_{i,j} = (1 - \delta_{i,j})A_{i,j}S_{i,j}.$$

Let us denote by $\widetilde{S}$ the matrix with off-diagonal entries $S_{i,j}$ and 0 on the diagonal. If we call $S = XX'/p$, we have

$$\widetilde{S} = S - \text{diag}(S).$$

Now note that $A_{i,j} = \psi_i + \psi_j$. Therefore, we have, if $\text{diag}(\psi)$ is the diagonal matrix with $(i,i)$ entry $\psi_i$,

$$H = \widetilde{S}\,\text{diag}(\psi) + \text{diag}(\psi)\widetilde{S}.$$

We just saw that under our assumptions, $\max_i |\psi_i| \to 0$ a.s. Because for any $n \times n$ matrices $L_1$, $L_2$, $\|\!\| L_1 L_2 \|\!\|_2 \leq \|\!\| L_1 \|\!\|_2 \|\!\| L_2 \|\!\|_2$, we see that to show that $\|\!\| H \|\!\|_2$ goes to 0, we just need to show that $\|\!\| \widetilde{S} \|\!\|_2$ remains bounded.

Now we clearly have, $\|\!\| S \|\!\|_2 \leq \|\!\| \Sigma \|\!\|_2 \|\!\| Y'Y/p \|\!\|_2$. We know from [45] that $\|\!\| Y'Y/p \|\!\|_2 \to \sigma^2(1 + \sqrt{n/p})^2$, a.s. Under our assumptions on $n$ and $p$, this is bounded. Now

$$\text{diag}(S) = \text{diag}(\psi) + \frac{\text{trace}(\Sigma)}{p}\,\text{Id},$$

so our concentration results once again imply that $\|\!\| \text{diag}(S) \|\!\|_2 \leq \text{trace}(\Sigma)/p + \eta$ a.s., for any $\eta > 0$. Because $\|\!\| \cdot \|\!\|_2$ is subadditive, we finally conclude that

$$\|\!\| \widetilde{S} \|\!\|_2 \text{ is bounded} \qquad \text{a.s.}$$

Therefore,

$$\|\!\| H \|\!\|_2 \to 0 \qquad \text{a.s.}$$

Putting together all these results, we see that we have shown that

$$\left\|\!\left\| T - (A \circ A - \text{diag}(A \circ A)) - 4\frac{\text{trace}(\Sigma^2)}{p^2}(11' - \text{Id}) \right\|\!\right\|_2 \to 0 \quad \text{in probability.}$$

• *Control of the third order term.* The third order term is the matrix $L$ with 0 on the diagonal and off-diagonal entries

$$L_{i,j} = \frac{f^{(3)}(\xi_{i,j})}{6}\left(\frac{\|X_i - X_j\|_2^2 - 2\,\text{trace}(\Sigma)}{p}\right)^3 \triangleq E \circ T,$$

where $T$ was the matrix investigated in the control of the second order term. On the other hand, $E$ is the matrix with entries

$$E_{i,j} = (1 - \delta_{i,j})\frac{f^{(3)}(\xi_{i,j})}{6}\left(\frac{\|X_i - X_j\|_2^2 - 2\operatorname{trace}(\Sigma)}{p}\right).$$

We have already seen that through concentration, we have

$$\max_{i \neq j}\left|\frac{\|X_i - X_j\|_2^2}{p} - \frac{2\operatorname{trace}(\Sigma)}{p}\right| \leq \log(p)p^{-1/10-\delta} \quad \text{a.s.}$$

This naturally implies that

$$\max_{i \neq j}\left|\xi_{i,j} - \frac{2\operatorname{trace}(\Sigma)}{p}\right| \leq \log(p)p^{-1/10-\delta} \quad \text{a.s.}$$

So if $f^{(3)}$ is bounded in a neighborhood of $\tau$, we see that with high-probability so is $\max_{i \neq j}|f^{(3)}(\xi_{i,j})|$. Therefore,

$$\max_{i \neq j}|E_{i,j}| \leq K\log(p)p^{-1/10-\delta}.$$

We are now in position to apply the Hadamard product argument (see Lemma A.5) we used for the control of the third order term in the proof of Theorem 2.1. To show that the third order term tends in operator norm to 0, we hence just need to show that $\|T\|_2$ remains small compared to the bound we just gave on $\max_{i,j}|E_{i,j}|$. Of course, this is equivalent to showing that the matrix that approximates $T$ has the same property in operator norm.

Clearly, because $\sigma_1(\Sigma)$ stays bounded, $\operatorname{trace}(\Sigma^2)/p$ stays bounded and so does $\|\operatorname{trace}(\Sigma^2)/p^2(11' - \operatorname{Id})\|_2$. So we just have to focus on $A \circ A - \operatorname{diag}(A \circ A)$. Recall that $A_{i,i} = 2(\|X_i\|_2^2/p - \operatorname{trace}(\Sigma)/p)$, and so $A_{i,i} = 2\psi_i$. We have already seen that our concentration arguments imply that $\max_i|\psi_i| \to 0$ a.s. So $\|\operatorname{diag}(A \circ A)\|_2 = \max_i \psi_i^2$ goes to 0 a.s. Now,

$$A = 1\psi' + \psi 1',$$

and hence, elementary Hadamard product computations [relying on $ab' \circ uv' = (a \circ u)(b \circ v)'$] give

$$A \circ A = 1(\psi \circ \psi)' + 2\psi\psi' + (\psi \circ \psi)1'.$$

Therefore,

$$\|A \circ A\|_2 \leq 2(\sqrt{n}\|\psi \circ \psi\|_2 + \|\psi\|_2^2).$$

Using Lemma A.1, and in particular equation (5), we see that

$$\mathbf{E}(\psi_i^2) = 2\sigma^4 \frac{\operatorname{trace}(\Sigma^2)}{p^2} + (\mu_4 - 3\sigma^4)\frac{\operatorname{trace}(\Sigma \circ \Sigma)}{p^2},$$

and therefore, $\mathbf{E}(\|\psi\|_2^2)$ remains bounded. On the other hand, using Lemma 2.7 of [5], we see that if we have $5 + \varepsilon$ absolute moments,

$$\mathbf{E}(\psi_i^4) \leq C\left(\frac{(\mu_4 \operatorname{trace}(\Sigma^2))^2}{p^4} + \mu_{5+\varepsilon} B_p^{3-\varepsilon} \frac{\operatorname{trace}(\Sigma^4)}{p^4}\right).$$

Now recall that we can take $B_p = p^{2/5-\delta}$. Therefore $n\mathbf{E}(\|\psi \circ \psi\|_2^2)$ is, at most, of order $B_p^{3-\varepsilon}/p$. We conclude that

$$P(\|\|A \circ A\|\|_2 > \log(p)\sqrt{B_p^{3-\varepsilon}/p}) \to 0.$$

Note that this implies that

$$P(\|\|T\|\|_2 > \log(p)\sqrt{B_p^{3-\varepsilon}/p}) \to 0.$$

Now, note that the third order term is of the form $E \circ T$. Because we have assumed that we have $5 + \varepsilon$ absolute moments, we have already seen that our concentration results imply that

$$\max_{i \neq j} |E_{i,j}| = \mathrm{O}\left(\log(p)\sqrt{\frac{B_p^2}{p}}\right) = \mathrm{O}(\log(p) p^{-1/10-\delta}) \qquad \text{a.s.}$$

Using the fact that $T$ has positive entries and therefore (see Lemma A.5), $\|\|E \circ T\|\|_2 \leq \max_{i,j} |E_{i,j}| \|\|T\|\|_2$, we conclude that with high-probability,

$$\|\|E \circ T\|\|_2 = \mathrm{O}\left((\log(p))^2 \sqrt{\frac{B_p^{5-\varepsilon}}{p^2}}\right) = \mathrm{O}((\log(p))^2 p^{-\delta'}) \qquad \text{where } \delta' > 0.$$

Hence,

$$\|\|E \circ T\|\|_2 \to 0 \qquad \text{in probability.}$$

• *Adjustment of the diagonal.* To obtain the compact form of the approximation announced in the theorem, we need to include diagonal terms that are not present in the matrices resulting from the Taylor expansion. Here, we show that the corresponding matrices are easily controlled in operator norm.

When we replace the zeroth and first order terms by

$$f(\tau)\mathbf{1}\mathbf{1}' + f'(\tau)\left[\mathbf{1}\psi' + \psi\mathbf{1}' - 2\frac{XX'}{p}\right],$$

we add to the diagonal the term $f(\tau) + f'(\tau)(2\psi_i - 2\|X_i\|_2^2/p) = f(\tau) - 2f'(\tau)\frac{\operatorname{trace}(\Sigma)}{p}$. In the end, we need to subtract it.

When we replace the second order term by $\frac{1}{2}f''(\tau)[1(\psi \circ \psi)' + 2\psi\psi' + (\psi \circ \psi)1' + 4\frac{\text{trace}(\Sigma^2)}{p^2}11']$, we add to the diagonal the diagonal matrix with $(i,i)$ entry,

$$2f''(\tau)\left[\psi_i^2 + \frac{\text{trace}(\Sigma^2)}{p^2}\right].$$

With our assumptions, $\max_i |\psi_i| \to 0$ a.s. and $\frac{\text{trace}(\Sigma^2)}{p}$ remains bounded, so the added diagonal matrix has operator norm converging to 0 a.s. We conclude that we do not need to add it to the correction in the diagonal of the matrix approximating our kernel matrix.  □

An interpretation of the proofs of Theorems 2.1 and 2.2 is that they rely on a local "multiscale" approximation of the original matrix (i.e., the terms used in the entry-wise approximation are all of different order of magnitudes, or at different "scales"). However, globally, that is, when looking at eigenvalues of the matrices and not just at each of their entries there is a bit of a mixture between the scales which creates the difficulties we had to deal with to control the second order term.

2.3.1. *A note on the Gaussian Kernel.* The Gaussian kernel corresponds to $f(x) = \exp(-\gamma x)$ in the notation of Theorem 2.2. We would like to discuss it a bit more because of its widespread use in applications.

The result of Theorem 2.2 gives accurate limiting eigenvalue information for the case where we renormalize the distances by the dimension which seems to be implicitly or explicitly what is often done in practice.

However, it is possible that information about the nonrenormalized case might also be of interest in some situations. Let us assume now that $\text{trace}(\Sigma)$ grows to infinity at least as fast as $p^{1/2+2/m+\delta}$ where $\delta > 0$ is such that $1/2 + 2/m + \delta < 1$ which is possible since $m \geq 5 + \varepsilon$ here. We of course still assume that its largest singular value, $\sigma_1(\Sigma)$, remains bounded. Then Corollary A.2 guarantees that

$$\min_{i \neq j} \frac{\|X_i - X_j\|_2^2}{p} > \frac{\text{trace}(\Sigma)}{p} \quad \text{a.s.}$$

Hence

$$\max_{i \neq j} \exp(-\|X_i - X_j\|_2^2) \leq \exp(-\text{trace}(\Sigma)) \leq \exp(-p^{1/2+2/m+\delta}) \quad \text{a.s.}$$

Hence, in this case, if $M$ is our kernel matrix with entries $\exp(-\|X_i - X_j\|_2^2)$, we have

$$\|\|M - \text{Id}\|\|_2 \leq n \exp(-p^{1/2+2/m+\delta}) \quad \text{a.s.,}$$

and the upper bound tends to zero extremely fast.

2.4. *More general models.* In this subsection, we consider more general models than the ones considered above. In particular, we will here focus on data models for which the vectors $X_i$ satisfy a so-called dimension-free concentration inequality. As was shown in [19], under these conditions, the Marčenko–Pastur equation holds (as well as generalized versions of it). Note that these models are more general than the one considered above (the proofs in the Appendix illustrate why the standard random matrix models can be considered as subcases of this more general class of matrices) and can describe various interesting objects like vectors with certain log-concave distributions or vectors sampled in a uniform manner from certain Riemannian submanifolds of $\mathbb{R}^p$ endowed with the canonical Riemannian metric inherited from $\mathbb{R}^p$.

Before we state more precisely the theorem and give examples of distributions that satisfy its assumptions, let us give some motivation. One potential criticism of Theorems 2.1 and 2.2 is that they deal with data models that are inherently quite linear. So a natural question is to understand whether our linear approximation result is limited to these linear settings. Also, kernel methods are often advocated for their handling of nonlinearities, and though the linear case is probably a basic one that needs to be understood, a null model of sorts, it is important to be able to go beyond it. As we will soon see, the next theorems allow us to get results beyond the linear setting.

Our generalization of Theorem 2.1 to more general distributions is the following.

THEOREM 2.3. *Consider a triangular "array" of matrices where each row of the array consists of a $n \times p$ matrix. We assume these matrices are independent. We call $X_i$, $i = 1, \ldots, n$ the rows of this matrix.*

*Suppose the vectors $\{X_i\}_{i=1}^n \in \mathbb{R}^p$ are i.i.d. mean 0 and have the property that for any 1-Lipschitz function $F$ (with respect to Euclidean norm), if $m_F$ is a median of $F(X_i)$,*

$$\forall t > 0 \qquad P(|F(X_i) - m_F| > t) \leq C \exp(-ct^b) \qquad \text{for some } b > 0,$$

*where $C$ is independent of $p$ and $c$ may depend on $p$ but is required to satisfy $c \geq p^{-(1/2-\varepsilon)b/2}$. $b$ is fixed and independent of $n$ and $p$.*

*Consider the $n \times n$ kernel random matrix $M$ with $M_{i,j} = f(X_i'X_j/p)$. Assume that $p \asymp n$.*

(a) *Call $\Sigma$ the covariance matrix of the $X_i$'s and assume that $\sigma_1(\Sigma)$ stays bounded and $\operatorname{trace}(\Sigma)/p$ has a limit.*
(b) *Suppose that $f$ is a real valued function which is $C^2$ around 0 and $C^1$ around $\operatorname{trace}(\Sigma)/p$.*

*Then the spectrum of this matrix is asymptotically nonrandom and has, a.s., the same limiting spectral distribution as that of*

$$\widetilde{M} = f(0)11' + f'(0)\frac{XX'}{p} + \upsilon_p \operatorname{Id}_n,$$

*where* $\upsilon_p = f(\frac{\operatorname{trace}(\Sigma)}{p}) - f(0) - f'(0)\frac{\operatorname{trace}(\Sigma)}{p}.$

We note that the term $f(0)11'$ does not affect the limiting spectral distribution of $\widetilde{M}$ since finite rank perturbations do not have any effect on limiting spectral distributions (see, e.g., [3], Lemma 2.2). Therefore, it could be removed from the approximating matrix, but since it will clearly be present in numerical work and simulations, we chose to leave it in our approximation. We also note that the limiting distribution of $XX'/p$ under these assumptions has been obtained in [19].

Here are a few examples of models satisfying the distributional assumptions stated above. (Unless otherwise noted, $b=2$ in all these examples.)

*Examples of distributions for which the previous theorem applies.*

- Gaussian random variables with $\|\!|\Sigma|\!\|_2$ bounded and $\operatorname{trace}(\Sigma)/p$ converges. The assumptions of the theorem apply according to [29], Theorem 2.7, with $c(p) = 1/\|\!|\Sigma|\!\|_2$.
- Vectors of the type $\sqrt{p}r$ where $r$ is uniformly distributed on the unit ($\ell_2$-) sphere is dimension $p$. Theorem 2.3 in [29] shows that Theorem 2.3 applies, with $c(p) = (1-1/p)/2$, after noticing that a 1-Lipschitz function with respect to Euclidean norm is also 1-Lipschitz with respect to the geodesic distance on the sphere.
- Vectors $\Gamma\sqrt{p}r$ with $r$ uniformly distributed on the unit ($\ell_2$-)sphere in $\mathbb{R}^p$ and with $\Gamma\Gamma' = \Sigma$ where $\Sigma$ satisfies the assumptions of the theorem.
- Vectors with log-concave density of the type $e^{-U(x)}$ with the Hessian of $U$ satisfying, for all $x$, $\operatorname{Hess}(U) \geq c\operatorname{Id}_p$ where $c > 0$ has the characteristics of $c(p)$ above (see [29], Theorem 2.7). Here we also need $\|\!|\Sigma|\!\|_2$ to satisfy the assumptions of the theorem.
- Vectors $r$ distributed according to a (centered) Gaussian copula, with corresponding correlation matrix $\Sigma$ having $\|\!|\Sigma|\!\|_2$ bounded. Here Theorem 2.3 applies, since if $\tilde{r}$ has a Gaussian copula distribution, then its $i$th entry satisfies $\tilde{r}_i = \Phi(v_i)$ where $v$ is multivariate normal with covariance matrix $\Sigma$, $\Sigma$ being a correlation matrix, that is, its diagonal is 1. Here $\Phi$ is the cumulative distribution function of a standard normal distribution which is trivially Lipschitz. Now taking $r = \tilde{r} - 1/2$ gives a centered Gaussian copula. The fact that the covariance matrix of $r$ then has bounded operator norm requires a bit of work and is shown in the Appendix of [19].

- Vectors sampled uniformly from certain compact connected smooth Riemannian submanifolds, $M$, of $\mathbb{R}^p$, canonically equipped with the Riemannian metric $g$ defined by restricting to each tangent space the ambient scalar product in $\mathbb{R}^p$. The curvature properties of these submanifolds need to satisfy the assumptions of Theorem 2.4 in [29]. Also, the covariance matrix $\Sigma$ need to satisfy the assumptions of Theorem 2.3. We note that since the length of a curve in $M$ is equal to its length in $\mathbb{R}^p$, the same remark that we made in the case of the sphere of unit radius applies here, too. In particular, a 1-Lipschitz function with respect to Euclidean norm is 1-Lipschitz with respect to the geodesic distance on the manifold.
- Vectors of the type $p^{1/b}r$, $1 \leq b \leq 2$ where $r$ is uniformly distributed in the 1-$\ell^b$ ball or sphere in $\mathbb{R}^p$. (See [29], Theorem 4.21, which refers to [34] as the source of the theorem.) We also refer the reader to [29], pages 37 and 38 for some of subtleties involved in the definition of the uniform distribution on the sphere. Fact A.1 applies to them, with $c(p)$ depending only on $b$. Also, the concentration function is of the form $\exp(-c(b)t^b)$ here.

We now turn to the proof of Theorem 2.3. The first step in the proof is the following lemma.

LEMMA 2.1. *Suppose $K_n$ is an $n \times n$ real symmetric matrix, whose spectral distribution converges weakly (i.e., in distribution) to a limit. Suppose $M_n$ is an $n \times n$ real symmetric matrix.*

1. *Suppose $M_n$ is such that $\|M_n - K_n\|_F = \mathrm{o}(\sqrt{n})$. Then $M_n$ and $K_n$ have the same limiting spectral distribution.*
2. *Suppose $M_n$ is such that $\|\|M_n - K_n\|\|_2 \to 0$. Then $M_n$ and $K_n$ have the same limiting spectral distribution.*

Before we prove the lemma, we note that our assumptions imply that the limiting spectral distribution of $K_n$ is a probability distribution. Therefore, to obtain the results of the lemma, we just need to show pointwise convergence of Stieltjes transforms and then rely on the results of [22], and in particular Corollary 1 there.

PROOF OF LEMMA 2.1. We call $\mathrm{St}_{K_n}$ and $\mathrm{St}_{M_n}$ the Stieltjes transforms of the spectral distributions of these two matrices. Suppose $z = u + iv$. Let us call $l_i(M_n)$ the $i$th largest eigenvalue of $M_n$.

*Proof of statement 1.* We first focus on the Frobenius norm part of the lemma. We have

$$|\mathrm{St}_{K_n}(z) - \mathrm{St}_{M_n}(z)| = \frac{1}{n}\left|\sum_{i=1}^{n} \frac{1}{l_i(K_n) - z} - \frac{1}{l_i(M_n) - z}\right| \leq \frac{1}{n}\sum_{i=1}^{n} \frac{|l_i(M_n) - l_i(K_n)|}{v^2}.$$

Now, by Holder's inequality,

$$\sum_{i=1}^{n} |l_i(M_n) - l_i(K_n)| \leq \sqrt{n} \sqrt{\sum_{i=1}^{n} |l_i(M_n) - l_i(K_n)|^2}.$$

Using Lidskii's theorem [i.e., the fact that, since $M_n$ and $K_n$ are hermitian, the vector with entries $l_i(M_n) - l_i(K_n)$ is majorized by the vector $l_i(M_n - K_n)$], with, in the notation of [10], Theorem III.4.4 $\Phi(x) = x^2$, we have

$$\sum_{i=1}^{n} |l_i(M_n) - l_i(K_n)|^2 \leq \sum_{i=1}^{n} l_i^2(M_n - K_n) = \|M_n - K_n\|_F^2.$$

We conclude that

$$|\mathrm{St}_{K_n}(z) - \mathrm{St}_{M_n}(z)| \leq \frac{\|M_n - F_n\|_F}{\sqrt{n}v^2},$$

since $|l_i(K_n) - z| \geq |\mathrm{Im}[l_i(K_n) - z]| = v$, and therefore $1/|l_i(K_n) - z| \leq 1/v$. Under the assumptions of the lemma, we therefore have

$$|\mathrm{St}_{K_n}(z) - \mathrm{St}_{M_n}(z)| \to 0.$$

Therefore the Stieltjes transform of the spectral distribution of $M_n$ converges pointwise to the Stieltjes transform of the limiting spectral distribution of $K_n$. Hence, by, e.g., Corollary 1 in [22], the spectral distribution of $M_n$ converges in distribution to the limiting spectral distribution of $K_n$, which, as noted earlier, is a probability distribution by our assumptions.

*Proof of statement 2.* Let us now turn to the operator norm part of the lemma. By the same computations as above, we have, using Weyl's inequality,

$$|\mathrm{St}_{K_n}(z) - \mathrm{St}_{M_n}(z)| = \frac{1}{n} \left| \sum_{i=1}^{n} \frac{1}{l_i(K_n) - z} - \frac{1}{l_i(M_n) - z} \right|$$

$$\leq \frac{1}{n} \sum_{i=1}^{n} \frac{|l_i(M_n) - l_i(K_n)|}{v^2}$$

$$\leq \frac{\|\|M_n - K_n\|\|_2}{v^2}.$$

Hence if $\|\|M_n - K_n\|\|_2 \to 0$, it is clear that the two Stieltjes transforms are asymptotically equal, and the conclusion follows. $\square$

We now turn to the proof of the theorem.

PROOF OF THEOREM 2.3. For the weaker statement required for the proof of Theorem 2.3, we will show that in the $\delta$-method we need to keep only the first term of the expansion as long as $f$ has a second derivative that is bounded in a neighborhood of 0, and a first derivative that is bounded in a neighborhood of $\text{trace}(\Sigma)/p$. In other words, we will split the problem into two parts: off the diagonal, we write

$$f\left(\frac{X_i'X_j}{p}\right) = f(0) + f'(0)\frac{X_i'X_j}{p} + \frac{f''(\xi_{i,j})}{2}\left(\frac{X_i'X_j}{p}\right)^2 \qquad \text{if } i \neq j;$$

on the diagonal, we write

$$f\left(\frac{X_i'X_i}{p}\right) = f\left(\frac{\text{trace}(\Sigma)}{p}\right) + f'(\xi_{i,i})\left(\frac{X_i'X_i}{p} - \frac{\text{trace}(\Sigma)}{p}\right).$$

• *Control of the off-diagonal error matrix.* Here we focus on the matrix $\widetilde{W}$ with $(i,j)$ entry

$$\widetilde{W}_{i,j} = 1_{i\neq j}\frac{f''(\xi_{i,j})}{2}\left(\frac{X_i'X_j}{p}\right)^2.$$

The strategy is going to be to control the Frobenius norm of the matrix

$$W_{i,j} = \begin{cases} \left(\dfrac{X_i'X_j}{p}\right)^2, & \text{if } i \neq j, \\ 0, & \text{if } i = j. \end{cases}$$

According to Lemma 2.1, it is enough for our needs to show that the Frobenius norm of this matrix is $o(\sqrt{n})$ a.s. to have the result we wish. Hence, the result will be shown, if we can for instance show that

$$\max_{i,j} W_{i,j} \leq p^{-(1/2+\varepsilon)}(\log(p))^{1+\delta} \qquad \text{a.s., for some } \delta > 0.$$

Now Lemma A.4 or Fact A.1 gives, for instance,

$$\max_{i\neq j}\left|\frac{X_i'X_j}{p}\right| \leq (pc^{2/b}(p))^{-1/2}[\log(p)]^{2/b} \qquad \text{a.s.}$$

Therefore, with our assumption on $c(p)$, we have

$$\max_{i,j} W_{i,j} \leq p^{-(1/2+\varepsilon)}(\log(p))^{4/b} \qquad \text{a.s.}$$

Now, $\|W\|_F \leq n\max_{i,j}|W_{i,j}|$, so we conclude that in this situation, with our assumptions that $n \asymp p$,

$$\|W\|_F = o(\sqrt{n}) \qquad \text{a.s.}$$

Now let us focus on

$$\widetilde{W}_{i,j} = f''(\xi_{i,j})W_{i,j},$$

where $\xi_{i,j}$ is between 0 and $X_i'X_j/p$. We just saw that with very high-probability, this latter quantity was less (in absolute value) than $p^{-(1/4+\varepsilon/2)} \times (\log(p))^{2/b}$, if $c \geq p^{-(1/2-\varepsilon)b/2}$. Therefore if $f''$ is bounded by $K$ in a neighborhood of 0, we have, with very high probability that

$$\|\widetilde{W}\|_F \leq K\|W\|_F = \mathrm{o}(\sqrt{n}).$$

• *Control of the diagonal matrix.* We first note that when we replace the off-diagonal matrix by $f(0)11' + f'(0)XX'/p$, we add to the diagonal certain terms that we need to subtract eventually.

Hence, our strategy here is to show that we can approximate (in operator norm) the diagonal matrix $D$ with entries

$$D_{i,i} = f\left(\frac{\mathrm{trace}(\Sigma)}{p}\right) + f'(\xi_{i,i})\left(\frac{X_i'X_i}{p} - \frac{\mathrm{trace}(\Sigma)}{p}\right) - f'(0)\frac{X_i'X_i}{p} - f(0)$$

by $v_p \mathrm{Id}_p$. To do so, we just have to show that the diagonal error matrix $Z$, with entries

$$Z_{i,i} = (f'(\xi_{i,i}) - f'(0))\left(\frac{X_i'X_i}{p} - \frac{\mathrm{trace}(\Sigma)}{p}\right)$$

goes to zero in operator norm.

As seen in Lemma A.4 or Fact A.1, if $c \geq p^{-(1/2-\varepsilon)b/2}$, with very high-probability,

$$\max_i \left|\frac{X_i'X_i}{p} - \frac{\mathrm{trace}(\Sigma)}{p}\right| \leq p^{-(1/4+\varepsilon/2)}(\log(p))^{2/b}.$$

If $f'$ is continuous and hence bounded around $\frac{\mathrm{trace}(\Sigma)}{p}$, we therefore see that the operator (or spectral) norm of $Z$ satisfies with high-probability

$$\|\|Z\|\|_2 \leq Kp^{-(1/4+\varepsilon/2)}(\log(p))^{2/b}.$$

• *Final step.* We clearly have

$$\widetilde{M} - M = W + Z.$$

It is also clear that $\widetilde{M}$ has a limiting spectral distribution, satisfying, up to centering and scaling, the Marčenko–Pastur equation; this was shown in [19]. By Lemma 2.1, we see that $\widetilde{M}$ and $\widetilde{M} - Z$ have the same limiting spectral distribution, since their difference is $Z$ and $\|\|Z\|\|_2 \to 0$. Using the same lemma, we see that $M$ and $\widetilde{M} - Z$ have (in probability) the same limiting spectral distribution, since their difference is $W$ and we have established that the Frobenius norm of this matrix is (in probability) $\mathrm{o}(\sqrt{n})$. Hence, $M$ and $\widetilde{M}$ have (in probability) the same limiting spectral distribution. $\square$

We finally treat the case of kernel matrices computed from Euclidean norms, in this more general distributional setting.

THEOREM 2.4. *Let us call $\tau = 2\operatorname{trace}(\Sigma)/p$ where $\Sigma$ is the covariance matrix of the $X_i$'s. Suppose that $f$ is a real valued function which is $C^2$ around $\tau$ and $C^1$ around 0.*

*Under the assumptions of Theorem 2.3, the kernel matrix $M$ with $(i,j)$ entry*

$$M_{i,j} = f\left(\frac{\|X_i - X_j\|_2^2}{p}\right)$$

*has a nonrandom limiting spectral distribution which is the same as that of the matrix*

$$\widetilde{M} = f(\tau)11' - 2f'(\tau)\frac{XX'}{p} + \upsilon_p \operatorname{Id}_n,$$

*where $\upsilon_p = f(0) + \tau f'(\tau) - f(\tau)$.*

We note once again that the term $f(\tau)11'$ does not affect the limiting spectral distribution of $M$. But we keep it for the same reasons as before.

PROOF OF THEOREM 2.4. Note that the diagonal term is simply $f(0)\operatorname{Id}$, so this term does not create any problem.

The rest of proof is similar to that of Theorem 2.3. In particular the control of the Frobenius norm of the second order term is done in the same way, by controlling the maximum of the off-diagonal term, using Corollary A.3 and Fact A.1 (and hence Lemma A.4).

Therefore, we only need to understand the first order term, in other words, the matrix with 0 on the diagonal and off-diagonal entry

$$\begin{aligned} R_{i,j} &= \frac{\|X_i - X_j\|_2^2}{p} - \tau \\ &= \left[\frac{\|X_i\|_2^2}{p} - \frac{\operatorname{trace}(\Sigma)}{p}\right] + \left[\frac{\|X_j\|_2^2}{p} - \frac{\operatorname{trace}(\Sigma)}{p}\right] - 2\frac{X_i' X_j}{p}. \end{aligned}$$

As in the proof of Theorem 2.2, let us call $\psi$ the vector with $i$th entry $\psi_i = \frac{\|X_i\|_2^2}{p} - \frac{\operatorname{trace}(\Sigma)}{p}$. Clearly,

$$R_{i,j} = \delta_{i,j}\left(1\psi' + \psi 1' - 2\frac{XX'}{p}\right).$$

Simple computations show that

$$R - 2\frac{\operatorname{trace}(\Sigma)}{p}\operatorname{Id} = 1\psi' + \psi 1' - 2\frac{XX'}{p}.$$

Now, obviously, $1\psi' + \psi 1'$ is a matrix of rank at most 2. Hence, $R$ has the same limiting spectral distribution as

$$2\frac{\text{trace}(\Sigma)}{p}\text{Id} - 2\frac{XX'}{p}$$

since finite rank perturbations do not affect limiting spectral distributions (see, for instance, [3], Lemma 2.2). This completes the proof. $\square$

2.5. *Some consequences of the theorems.* In practice, it is often the case that slight variant of kernel random matrices are used. In particular, it is customary to center the matrices, that is, to transform $M$ so that its row sum, or column sum or both are 0. Note that these operations correspond to right and/or left multiplication by the matrix $H = \text{Id}_n - 11'/n$.

In these situations, our results still apply; the following fact makes it clear.

FACT 2.1 (Centered kernel random matrices). *Let $H$ be the $n \times n$ matrix $\text{Id}_n - 11'/n$.*

1. *If the kernel random matrix $M$ can be approximated consistently in operator norm by $K$, then, if $a, b \in \{0, 1\}$,*

   $H^a M H^b$ *can be approximated consistently in operator norm by $H^a K H^b$.*

2. *If the kernel random matrix $M$ has the same limiting spectral distribution as the matrix $K$, then, if $a, b \in \{0, 1\}$,*

   $H^a M H^b$ *has the same limiting spectral distribution as $K$.*

A nice consequence of the first point is that the recent hard work on localizing the largest eigenvalues of sample covariance matrices (see [8, 17] and [31]) can be transferred to kernel random matrices and used to give some information about the localization of the largest eigenvalues of $HMH$, for instance. In the case of the results of [17], Fact 2 and the arguments of [19], Section 2.3.4, show that it gives exact localization information. In other words, we can characterize the a.s. limit of the largest eigenvalue of $HMH$ (or $HM$ or $MH$) fairly explicitly, provided Fact 2 in [17] applies. Finally, let us mention the obvious fact that since two square matrices $A$ and $B$, $AB$ and $BA$ have the same eigenvalues, we see that $HMH$ has the same eigenvalues as $MH$ and $HM$ because $H^2 = H$.

PROOF OF FACT 2.1. The proofs are simple. First note that $H$ is positive semi-definite and $|||H|||_2 = 1$. Using the submultiplicativity of $||| \cdot |||_2$, we see that

$$|||H^a M H^b - H^a K H^b|||_2 \leq |||M - K|||_2 |||H^a|||_2 |||H^b|||_2 = |||M - K|||_2.$$

This shows the first point of the fact.

The second point follows from the fact that $H^a M H^b$ is a finite rank perturbation of $M$. Hence, using Lemma 2.2 in [3], we see that these two matrices have the same limiting spectral distribution, and since, by assumption, $K$ has the same limiting spectral distribution as $M$, we have the result of the second point. $\square$

*On Laplacian-like matrices.* Finally, we point out a simple consequence of our results for Laplacian-like matrices. In light of recent results on manifold learning (see [9]) where these matrices play a key role, it is natural to ask what happens to them in our context. Suppose $M$ is an $n \times n$ kernel random matrix as defined in the previous theorems, and consider the (Laplacian-like) matrix $L$ defined by

$$\begin{cases} L_{i,j} = -M_{i,j}/n, & \text{if } i \neq j, \\ L_{i,i} = \dfrac{1}{n} \sum_{i \neq j} M_{i,j}. \end{cases}$$

Call $D_L$ the diagonal matrix made up of the diagonal elements of $L$. We note that our concentration results (Lemmas A.3 and A.4, as well as Fact A.1) imply that $D_L$ can be approximated in operator norm by $f(0) \operatorname{Id}_n$ in the scalar product kernel matrix case and by $f(2 \operatorname{trace}(\Sigma_p)/p) \operatorname{Id}_n$ in the Euclidean distance kernel matrix case. [It is so because $L_{i,i}$ is an average of almost constant (and equal) quantities, so with high-probability $L_{i,i}$ cannot deviate from this constant value, for that would require that at least one of the components of the average deviate from the constant value in question.] Hence, there exists $\gamma_p$ such that $\|\!|D_L - \gamma_p \operatorname{Id}_n|\!\|_2$ tends to 0 almost surely. We also recall that the diagonal of the matrix $M$ can be consistently approximated in operator norm by a (finite) multiple of the identity matrix, so the diagonal of $M/n$ can be consistently approximated in operator norm by 0. Therefore, $\|\!|L + M/n - D_L|\!\|_2$ tends to 0 almost surely, and therefore, $\|\!|L + M/n - \gamma_p \operatorname{Id}_n|\!\|_2$ tends to zero almost surely. In other words, $L$ can be consistently approximated in operator norm by $\gamma_p \operatorname{Id}_n - M/n$. Consequently, when we can, as in Theorems 2.1 and 2.2, consistently approximate $M$ in operator norm by a linearized version, $K$, of $M$, then $L$ can be consistently approximated in operator norm by $\gamma_p \operatorname{Id}_n - K/n$, and we can deduce spectral properties of $L$ from that of $\gamma_p \operatorname{Id}_n - K/n$. When we know only about the limiting spectral distribution of $M$, as in Theorems 2.3 and 2.4, the operator norm consistent approximation of $L$ by $\gamma_p \operatorname{Id}_n - M/n$ carries over to give us information about the limiting spectral distribution of $L$ since the effect of $\gamma_p \operatorname{Id}_n$ is just to "shift" the eigenvalues by $\gamma_p$. We note that getting information about the eigenvectors of $L$ would require finer work on the properties of the matrix $D_L$ since approximating it by a multiple of the identity does not give us any information about its eigenvectors.

**3. Conclusions.** The main result of this paper is that under various technical assumptions, in high-dimensions, kernel random matrices [i.e., $n \times n$ matrices with $(i,j)$th entry $f(X_i'X_j/p)$ or $f(\|X_i - X_j\|_2^2/p)$ where $\{X_i\}_{i=1}^n$ are i.i.d. random vectors in $\mathbb{R}^p$ with $p \to \infty$ and $p \asymp n$] which are often used to create nonlinear versions of standard statistical methods and essentially behave like covariance matrices, that is, linearly, a result that is in sharp contrast with the low-dimensional situation where $p$ is assumed to be fixed, and where it is known that, under some regularity conditions, spectral properties of kernel random matrices mimick those of certain integral operators. Under ICA-like assumptions, we were able to get a "strong approximation" result (Theorems 2.1 and 2.2), that is, an operator norm consistency result that carries information about individual eigenvalues and eigenvectors corresponding to separated eigenvalues. Under more general and less linear assumptions (Theorems 2.3 and 2.4), we have obtained results concerning the limiting spectral distribution of these matrices using a "weak approximation" result relying on bounds on Frobenius norms.

Beside the mathematical results obtained above, this study raises several statistical questions, both about the richness—or lack thereof—of models that are often studied in random matrix theory and about the effect of kernel methods in this context.

3.1. *On kernel random matrices.* Our study, motivated in part by numerical experiments we read about in the interesting [44], has shown that in the asymptotic setting we considered which is generally considered relevant for high-dimensional data analysis, the kernel random matrices considered here behave essentially like matrices closely connected to sample covariance matrices. This is in sharp contrast to the low-dimensional setting where it was explained heuristically in [44], and proved rigorously in [28], that the eigenvalues of kernel random matrices converged (under certain assumptions) to those of a canonically related operator.

This suggests that kernel methods could suffer from the same problems that affect linear statistical methods, such as Principal Component Analysis, in high-dimensions. The practical significance of our result is that in high-dimensions, the nonlinear methods that rely on kernel matrices may be behaving like their linear counterparts. Our study also permits the transfer of some recent random matrix results concerning large-dimensional sample covariance matrices to kernel random matrices. We now discuss some possible practical settings to highlight when our results are and are not relevant.

*On kernel-PCA.* An important motivation for this study was to try to understand the properties of kernel-PCA in high-dimensions. We refer the reader to [36] pages 48–50 for a primer on kernel-PCA, but let us say that kernel-PCA performs a spectral decomposition of a (row and column-centralized) kernel matrix to efficiently perform a nonlinear version of PCA; instead of doing standard PCA, the algorithm performs PCA in feature space. Our results, and in particular Theorems 2.1 and 2.2 clearly show that in high-dimensions, when the assumptions of the theorems are satisfied, the algorithm may essentially be performing a linear PCA despite appearances to the contrary. By contrast, in low-dimension, results such as [28] show that the intuition behind kernel-PCA is correct and that the algorithm then performs a genuinely nonlinear PCA. Being aware of the difference between the two settings should be helpful to practitioners in that it will inform them about the possible limitations of kernel-PCA as a nonlinear method. For an example of applications, we refer the reader to, for instance, [44] and [36], Chapter 10. We note also that from a slightly more "numerical analysis" standpoint, our results basically say that the Nyström method (see [44] and references therein) to approximate eigenfunctions of integral operators can be unreliable in high-dimensions.

*On kernel-ICA.* The setting of Theorems 2.1 and 2.2 is naturally well-suited for applications to ICA-like problems. Since we are dealing with vectors $X_i$ here, we would be considering multidimensional ICA problems (see [2], page 33). For instance, in the formulation of [2], kernel-ICA is solved by solving a kernel-CCA problem, that is, a generalized eigenvalue problem with kernel matrices as input. We refer the reader to equations (10) and (13) in [2] for more details. The results of Theorems 2.1 and 2.2 are directly relevant here, since the matrices at stake in these kernel-CCA problems can be approximated consistently in operator norm by linear counterparts, and hence the solution of the kernel-CCA problem can be consistently approximated by the solution of the problem obtained by linearizing the kernel matrices at stake, provided the smallest singular value of the linearized version of the kernel matrices in question stay bounded away from 0. We note that in practice this latter requirement can be checked using our fairly detailed knowledge of the properties of extreme eigenvalues of sample covariance matrices. We note that in this setting, our theorems confirm the predictions in [2], page 33, that problems might arise with the algorithm they propose in high-dimensions due to slow decay of eigenvalues.

*Geostatistics applications.* Certain kernel matrices, corresponding to covariance functions of, for instance, Gaussian processes, also appear in geostatistics and spatial statistics in techniques such as kriging (see, e.g., [15], Chapter 3 and, for instance, pages 106–110). For examples of kernels that correspond to covariance functions of Gaussian processes, we refer the reader to [33], Chapter 4. Naturally, in those sort of applications, the dimension of the data vectors is low (at most 3), and therefore "classical results" such as those of [28] apply whereas our results are limited to the high-dimensional setting more often encountered in some problems of multivariate statistics, machine learning or bioinformatics.

3.2. *Limitations of standard random matrix models.* In the study of spectral distribution of large-dimensional sample covariance matrices, it has been somewhat forcefully advocated that the study should be done under the assumptions that the data are of the form $X_i = \Sigma^{1/2} Y_i$ where the entries of $Y_i$ have, for instance, finite fourth moment. At first sight, this idea is appealing, as it seems to allow a great variety of distributions and hence flexible modeling. A possible drawback however, is the assumption that the data are linear combinations of i.i.d. random variables or the necessary presence of independence in the model. This has however been recently addressed (see, e.g., [19]) and it has been shown that one could go beyond models requiring independence in a lurking random vector which the data linearly depend on.

*Data analytic consequences.* However, a serious limitation is still present. As the results of Lemmas A.3, A.4 and Fact A.1 make clear, under the models for which the limiting spectral distribution of the sample covariance matrix has been shown to satisfy the Marčenko–Pastur equation, the norms of the data vectors are concentrated, and the corresponding data vectors are almost orthogonal to one another. In other words, under the "standard" ICA-like random matrix models (used in Theorems 1 and 2), that is, the random vectors $\{X_i\}_{i=1}^n$ are i.i.d. in $\mathbb{R}^p$, $X_i = \Sigma^{1/2} Y_i$ with $Y_i$ having i.i.d. entries with mean 0, variance 1 and $4 + \varepsilon$ absolute moments, we have, assuming that $\{Y_i\}_{i=1}^n$ are i.i.d., $p \to \infty$, $p/n$ and $\|\!|\Sigma|\!\|_2$ remain bounded; the vectors $\{X_i\}_{i=1}^n$ have the property that for a deterministic (and computable) sequence $c_p$, we have

$$\max_{1 \leq i \leq n} |\|X_i\|_2^2 / p - c_p| \to 0$$

and

$$\max_{i \neq j} |X_i' X_j| / p \to 0.$$

Both these statements hold almost surely. Geometrically, this means that the vectors $\{X_i / \sqrt{p}\}_{i=1}^n$ are close to a sphere and almost orthogonal to one another. These properties also hold for the more general (and less linear) models we considered in Theorems 2.3 and 2.4.

Hence, if one were to plot a histogram of $\{\|X_i\|_2^2 / p\}_{i=1}^n$, this histogram would look tightly concentrated around a single value—the spread of this histogram being computable from our concentration results (Lemmas A.3, A.4 and Fact A.1). Though the models appear to be quite rich, the geometry that we can perceive by sampling $n$ such vectors, with $n \asymp p$, is, arguably, relatively poor. These remarks should not be taken as aiming to discredit the interesting body of work that has emerged out of the study of such models. Their aim is just to warn possible users that in data analysis, a good first step would be to plot the histogram of $\{\|X_i\|_2^2/p\}_{i=1}^n$ and check whether it is concentrated around a single value. Similarly, one might want to plot the histogram of inner products $\{X_i'X_j/p\}$ and check that it is concentrated around 0. If this is not the case, then insights derived from random matrix theoretic studies would likely not be helpful in the data analysis.

We note, however, that recent random matrix work (see [13, 14, 19, 32]) has been concerned with distributions which could be loosely speaking be called of "elliptical" type—though they are more general than what is usually called elliptical distributions in Statistics. In those settings, the data is, for instance, of the form $X_i = r_i \Sigma^{1/2} Y_i$ where $r_i$ is a real-valued random variable, independent of $Y_i$. This allows the data vectors to not approximately live on spheres (but does not change anything about angles between different vectors), and is a possible way to address some of the concerns we just raised. The characterization of the limiting spectrum gets quite a bit more involved than in the "standard" setting, that is, $r_i = 1$, and the results show a lack of robustness to the "indirect" assumption that the data vectors live close to a sphere.

Finally, this geometric discussion applies also to theoretical studies undertaken under the assumptions that the $X_i$ are $\mathcal{N}(0, \Sigma_p)$ and that the problem is high dimensional. It should highlight some possibly severe limitations of the normality assumption in high-dimensions.

## APPENDIX

In this appendix, we collect a few useful results that are needed in the proof of our theorems, and whose content we thought would be more accessible if they were separated from the main proofs.

**(A) Some useful results.** We have the following elementary facts.

LEMMA A.1. *Suppose $Y$ is a vector with i.i.d. entries and mean 0. Call its entries $y_i$. Suppose $\mathbf{E}(y_i^2) = \sigma^2$ and $\mathbf{E}(y_i^4) = \mu_4$. Then if $M$ is a deterministic matrix,*

(4) $\mathbf{E}(YY'MYY') = \sigma^4(M + M') + (\mu_4 - 3\sigma^4)\operatorname{diag}(M) + \sigma^4 \operatorname{trace}(M)\operatorname{Id}.$

*Further, we have $(Y'MY)^2 = \operatorname{trace}(MYY'MYY')$ and*

(5) $\mathbf{E}(\operatorname{trace}(MYY'MYY'))$
$= \sigma^4 \operatorname{trace}(M^2 + MM') + \sigma^4(\operatorname{trace}(M))^2 + (\mu_4 - 3\sigma^4)\operatorname{trace}(M \circ M).$

*Here* $\mathrm{diag}(M)$ *denotes the matrix consisting of the diagonal of the matrix* $M$ *and* 0 *off the diagonal. The symbol* ∘ *denotes Hadamard multiplication between matrices.*

PROOF. Let us call $R = YY'MYY'$. The proof of the first part is elementary and consists merely in writing the $(i,j)$th entry of the corresponding matrix. As a matter of fact, we have

$$R_{i,j} = y_i y_j \sum_{i,j} y_i y_j M_{i,j} = \sum_{k,l} y_i y_j y_k y_l M_{k,l}.$$

Using the fact that entries of $Y$ are independent and have mean 0, we see that, in the sum, the only terms that will not be 0 in expectation are those for which each index appears at least twice. If $i \neq j$, only the terms of the form $y_i^2 y_j^2$ have this property. So if $i \neq j$,

$$\mathbf{E}(R_{i,j}) = \mathbf{E}(y_i^2 y_j^2 (M_{i,j} + M_{j,i})) = \sigma^4 (M_{i,j} + M_{j,i}).$$

Let us now turn to the diagonal terms. Here again, only the terms $y_i^2 y_k^2$ matter. So on the diagonal,

$$\mathbf{E}(R_{i,i}) = \mu_4 M_{i,i} + \sigma^4 \sum_{j \neq i} M_{j,j} = (\mu_4 - \sigma^4) M_{i,i} + \mathrm{trace}(M).$$

We conclude that

$$\mathbf{E}(R) = \sigma^4 (M + M') + (\mu_4 - 3\sigma^4) \mathrm{diag}(M) + \mathrm{trace}(M) \mathrm{Id}.$$

The second part of the proof follows from the first result, after we remark that, if $D$ is a diagonal and $L$ is general matrix, $\mathrm{trace}(LD) = \mathrm{trace}(L \circ D)$, from which we conclude that $\mathrm{trace}(M \mathrm{diag}(M)) = \mathrm{trace}(M \circ \mathrm{diag}(M)) = \mathrm{trace}(M \circ M)$. □

LEMMA A.2 (Concentration of quadratic forms). *Suppose the vectors* $Z$ *is a vector in* $\mathbb{R}^p$ *with i.i.d. entries of mean 0 and variance* $\sigma^2$. *Suppose that their entries are bounded by* $B_p$. *Let* $M$ *be a symmetric matrix, with largest singular value* $\sigma_1(M)$. *Call*

$$\zeta_p = \frac{128 \exp(4\pi) \sigma_1(M) B_p^2}{p},$$

$$\nu_p = \sqrt{\sigma_1(\Sigma)}.$$

*Then we have, if* $r/2 > \zeta_p$,

$$P\left( \left| \frac{Z'MZ}{p} - \sigma^2 \frac{\mathrm{trace}(M)}{p} \right| > r \right)$$

(6)
$$\leq 8 \exp(4\pi) \exp(-p(r/2 - \zeta_p)^2 / (32 B_p^2 (1 + 2\nu_p)^2 \sigma_1(M)))$$
$$+ 8 \exp(4\pi) \exp(-p/(32 B_p^2 (1 + 2\nu_p)^2 \sigma_1(M))).$$

PROOF. We can decompose, using the spectral decomposition of $M$, $M = M_+ - M_-$ where $M_+$ is positive semi-definite and $M_-$ is positive definite (or 0 if $M$ is itself positive semi-definite). We can do so by replacing the negative eigenvalues of $M$ by 0 in the spectral decomposition and get $M_+$ in that way. Note that then, the largest singular values of $M_+$ and $M_-$ are also bounded by $\sigma_1(M)$ since $\sigma_1(M)$ is absolute value of the largest eigenvalue of $M$ in absolute value, and the nonzero eigenvalues of $M_+$ are a subset of the eigenvalues of $M$ and so are the eigenvalues of $M_-$ when $M_-$ is not 0. Now it is clear that the function $F$ which associates to a vector $x$ in $\mathbb{R}^p$ the scalar $\sqrt{x'M_+x/p} = \|M_+^{1/2}x/\sqrt{p}\|_2$ is a convex, $\sqrt{\sigma_1(M)/p}$-Lipschitz function with respect to Euclidean norm. Calling $m_F$ the median of $F(Z)$, we have, using Corollary 4.10 in [29],

$$P(|F(Z) - m_F| > r) \leq 4\exp(-pr^2/(16B_p^2\sigma_1(M))).$$

Let us now call $\mu_F$ the mean of $F(Z)$ (it exists according to Proposition 1.8 in [29]). Following the arguments given in the proof of this Proposition 1.8, and spelling out the constants appearing in the last result of Proposition 1.8 in [29], we see that

$$P(|F(Z) - \mu_F| > r) \leq 4\exp(4\pi)\exp(-pr^2/(32B_p^2\sigma_1(M))).$$

(Using the notation of Proposition 1.8 in [29], we picked $\kappa_2 = 1/2$, and $C' = \exp(\pi C^2/4)$; showing that this is a valid choice just requires one to carry out some of the computations mentioned in the proof of that Proposition.)

Let us call $A, B, D$ the sets

$$A \triangleq \left\{ \left| \frac{Z'M_+Z}{p} - \mu_F^2 \right| > r \right\},$$

$$B \triangleq \left\{ \sqrt{\frac{Z'M_+Z}{p}} + \mu_F \leq 1 + 2\mu_F \right\} = \left\{ \sqrt{\frac{Z'M_+Z}{p}} - \mu_F \leq 1 \right\}$$

and

$$D \triangleq \left\{ \left| \sqrt{\frac{Z'M_+Z}{p}} - \mu_F \right| > r/(1 + 2\mu_F) \right\}.$$

Of course, we have $P(A) \leq P(A \cap B) + P(B^c)$. Now note that $A \cap B \subseteq D$, simply because for positive reals, $a - b/(\sqrt{a} + \sqrt{b}) = \sqrt{a} - \sqrt{b}$. We conclude that

$$P(A) \leq 4\exp(4\pi)[\exp(-pr^2/(32B_p^2(1+2\mu_F)^2\sigma_1(M)))$$
$$+ \exp(-p/(32B_p^2\sigma_1(M)))].$$

Let us know call $\sigma^2$ the variance of the each of the component of $Z$. We know, according to Proposition 1.9 in [29], that

$$\operatorname{var}(F) = \frac{\mathbf{E}(Z'M_+Z)}{p} - \mu_F^2 = \sigma^2 \frac{\operatorname{trace}(M_+)}{p} - \mu_F^2 \leq \zeta_p$$
$$= \frac{128\exp(4\pi)\sigma_1(M)B_p^2}{p}.$$

Hence, we conclude that, if $r > \zeta_p$,

$$P\left(\left|\frac{Z'M_+Z}{p} - \sigma^2\frac{\operatorname{trace}(M_+)}{p}\right| > r\right)$$
$$\leq 4\exp(4\pi)\exp(-p(r-\zeta_p)^2/(32B_p^2(1+2\mu_F)^2\sigma_1(M)))$$
$$+ 4\exp(4\pi)\exp(-p/(32B_p^2(1+2\mu_F)^2\sigma_1(M))).$$

To get the announced result, we note that for the sum of two reals to be greater than $r$ in absolute value, one needs to be greater than $r/2$, and that our bounds become conservative when we replace $\mu_F$ (and its counterpart for $M_-$) by $\nu_p$. [Note that we get conservative bounds when replacing the $\mu_F$'s by $\max(\mathbf{E}(\sqrt{Z'M_+Z/p}), \mathbf{E}(\sqrt{Z'M_-Z/p}))$, and that this quantity is clearly bounded by $\sigma\sigma_1(\Sigma)$.] Hence, we have, as announced, if $r/2 > \zeta_p$,

$$P\left(\left|\frac{Z'MZ}{p} - \sigma^2\frac{\operatorname{trace}(M)}{p}\right| > r\right)$$
$$\leq 8\exp(4\pi)\exp(-p(r/2-\zeta_p)^2/(32B_p^2(1+2\mu_F)^2\sigma_1(M)))$$
$$+ 8\exp(4\pi)\exp(-p/(32B_p^2(1+2\mu_F)^2\sigma_1(M))).$$

Finally, we note that the proof makes clear that the same result would hold for different choices of $M_+$ and $M_-$, as long as $\max(\sigma_1(M_+), \sigma_1(M_-)) \leq \sigma_1(M)$. □

We therefore have the following useful corollary:

COROLLARY A.1. *Let $Y_i$ and $Y_j$ be i.i.d. random vectors as in Lemma A.2 with variance 1. Suppose that $\Sigma$ is a positive semi-definite matrix. We have, with*

$$\zeta_p = \frac{128\exp(4\pi)\sigma_1(\Sigma)B_p^2}{p}$$

*and*

$$\nu_p = \sqrt{\sigma_1(\Sigma)},$$

*that if $r/2 > \zeta_p$ and $K = 8\exp(4\pi)$,*

(7)
$$P\left(\left|\frac{Y_i'\Sigma Y_j}{p}\right| > r\right) \leq K\exp(-p(r/2-\zeta_p)^2/(32B_p^2(1+2\nu_p)^2\sigma_1(\Sigma)))$$
$$+ K\exp(-p/(32B_p^2(1+2\nu_p)^2\sigma_1(\Sigma))).$$

PROOF. The proof relies on the results of Lemma A.2. Remark that, since $\Sigma$ is symmetric,

$$Y_i'\Sigma Y_j = \frac{1}{2}(Y_i'Y_j')\begin{pmatrix} 0 & \Sigma \\ \Sigma & 0 \end{pmatrix}\begin{pmatrix} Y_i \\ Y_j \end{pmatrix}.$$

Now the entries of the vector made by concatenating $Y_i$ and $Y_j$ are i.i.d. and so we fall back into the setting of Lemma A.2. Finally, here $M_+$ and $M_-$ are known explicitly. A possible choice is $M_+ = 1/2\begin{pmatrix} \Sigma & \Sigma \\ \Sigma & \Sigma \end{pmatrix}$ and $M_- = 1/2\begin{pmatrix} \Sigma & 0 \\ 0 & \Sigma \end{pmatrix}$. $\nu_p$ is obtained by upper bounding the expectation of the square of $F$ in the notation of the proof of the previous lemma for these explicit matrices. Note that their largest singular values are both smaller that $\sigma_1(\Sigma)$, so the results of the previous lemma apply. $\square$

LEMMA A.3. *Let $\{Y_i\}_{i=1}^n$ be i.i.d. random vectors in $\mathbb{R}^p$, whose entries are i.i.d., mean 0, variance 1 and have bounded (in p) $m \geq 4$ absolute moments. Suppose that $\{\Sigma_p\}$ is a sequence of positive semi-definite matrices whose operator norms are uniformly bounded in p and $n/p$ is asymptotically bounded. We have, for any given $\varepsilon > 0$,*

$$\max_{i,j}\left|\frac{Y_i'\Sigma_p Y_j}{p} - \delta_{i,j}\frac{\mathrm{trace}(\Sigma_p)}{p}\right| \leq p^{-1/2+2/m}(\log(p))^{(1+\varepsilon)/2} \qquad a.s.$$

PROOF. Throughout the proof, we assume without loss of generality that $m < \infty$.

Call $t = 2/m$. It is clear that with our moment assumptions, $t \leq 1/2$. According to Lemma 2.2 in [45], the maximum of the array of $\{Y_i\}_{i=1}^n$ is a.s. less than $p^t$. So to control the maximum of the inner products of interest, it is enough to control the same quantity when we replace $Y_i$ by $\widetilde{Y}_i$ with $\widetilde{Y}_{i,l} \triangleq Y_{i,l}1_{|Y_{i,l}|\leq p^t}$. Now note that $\widetilde{Y}_i$ satisfies the boundedness assumption of Corollary A.1, but its mean is not necessarily zero and its variance is not 1. Note however, that all the entries of $\widetilde{Y}_i$ have the same mean, $\widetilde{\mu}$. Since $Y_i$ has mean 0, we have

$$|\widetilde{\mu}| \leq \mathbf{E}(|Y_{1,1}|1_{|Y_{1,1}|>p^t}) \leq \mathbf{E}(|Y_{1,1}|^m p^{-t(m-1)}) \leq \mu_m p^{-2+t}.$$

Similarly, if we call $\widetilde{\sigma}^2$ the variance of $\widetilde{Y}$, we have

$$\widetilde{\sigma}^2 = \mathbf{E}(|Y_{1,1}|^2 1_{|Y_{1,1}|\leq p^t}) - \widetilde{\mu}^2 = 1 - (\mathbf{E}(|Y_{1,1}|^2 1_{|Y_{1,1}|>p^t}) + \widetilde{\mu}^2).$$

Hence, $0 \leq 1 - \widetilde{\sigma}^2$, and

$$\begin{aligned} 1 - \widetilde{\sigma}^2 &= \mathbf{E}(|Y_{1,1}|^2 1_{|Y_{1,1}| > p^t}) + \widetilde{\mu}^2 \\ &\leq \mathbf{E}(|Y_{1,1}|^m p^{-t(m-2)}) + \widetilde{\mu}^2 \\ &\leq \mu_m p^{-2+2t} + \mu_m^2 p^{-4+2t} = \mathrm{O}(p^{-2+2t}). \end{aligned}$$

Let us call $U_i = \widetilde{Y}_i - \widetilde{\mu} 1_p$ and $\widetilde{U}_i = U_i/\widetilde{\sigma}$ where $\widetilde{\sigma}^2$ is also the variance of $U_i$. Corollary A.1 applies to the random variables $\widetilde{U}_i$ with $B_p = 2p^t$ when $p$ is large enough. So $\zeta_p = \mathrm{O}(p^{1-2t})$. Let us now call, for some $\varepsilon > 0$,

$$r(p) = p^{t-1/2}(\log(p))^{(1+\varepsilon)/2}.$$

Since, for $p$ large enough, $r(p)/2 > \zeta_p$, we can apply the conclusions of Corollary A.1, and by plugging in the different quantities, we see that

$$P(|\widetilde{U}_i' \Sigma_p \widetilde{U}_j / p| > r(p)) \leq \exp(-K(\log(p))^{1+\varepsilon}),$$

where $K$ denotes a generic constant (that may change from display to display). In particular, $K$ is independent of $p$ and is hence trivially bounded away from 0 as $p$ grows. The bound we just obtained on $1 - \widetilde{\sigma}^2$ also implies that for $p$ large enough, $\widetilde{\sigma}^2 > 1/2$ from which we conclude that for another $K$ with the same properties,

$$P(|U_i' \Sigma_p U_j / p| > r(p)) \leq \exp(-K(\log(p))^{1+\varepsilon}).$$

In other respects, the arguments of Lemma A.2 show that, since $\widetilde{\sigma}^2$ is the variance of $U_i$,

$$P(|U_i' \Sigma_p U_i / p - \widetilde{\sigma}^2 \operatorname{trace}(\Sigma_p)/p| > r(p)) \leq \exp(-K(\log(p))^{1+\varepsilon}).$$

Now

$$\frac{\widetilde{Y}_i' \Sigma_p \widetilde{Y}_j}{p} = \frac{U_i' \Sigma_p U_j}{p} + \widetilde{\mu} \frac{(1' \Sigma_p U_j + U_i' \Sigma_p 1)}{p} + \widetilde{\mu}^2 \frac{1' \Sigma_p 1}{p}.$$

Remark that $1' \Sigma_p 1 \leq p \sigma_1(\Sigma_p)$, and $|1' \Sigma_p U_j| \leq \sqrt{1' \Sigma_p 1} \sqrt{U_j' \Sigma_p U_j}$. We conclude, using the results obtained in the proof of Lemma A.2 that with probability greater than $1 - \exp(-K(\log(p))^{1+\varepsilon})$, the middle term is smaller than $2\sqrt{\sigma_1(\Sigma_p)}(\sqrt{\sigma_1(\Sigma_p)} + r(p))|\widetilde{\mu}|$. As a matter of fact, $\sqrt{U_j' \Sigma_p U_j/p}$ is concentrated around its mean which is smaller than $\widetilde{\sigma}\sqrt{\operatorname{trace}(\Sigma_p)/p}$ which is itself smaller than $\sqrt{\sigma_1(\Sigma_p)}$. Now recall that $|\widetilde{\mu}| = \mathrm{O}(p^{-2+t}) = \mathrm{o}(r(p))$. We can therefore conclude that

$$P\left(\left|\frac{\widetilde{Y}_i' \Sigma_p \widetilde{Y}_j}{p} - \delta_{i,j} \widetilde{\sigma}^2 \frac{\operatorname{trace}(\Sigma_p)}{p}\right| > 2r(p)\right) \leq 2\exp(-K(\log(p))^{1+\varepsilon}).$$

Now note that $0 \leq 1 - \widetilde{\sigma}^2 = \mathrm{O}(p^{-2+2t}) = \mathrm{o}(r(p))$ since $t \leq 1/2 < 3/2$. With our assumptions, $\mathrm{trace}(\Sigma_p)/p$ remains bounded, so we have finally

$$P\left(\left|\frac{\widetilde{Y}_i' \Sigma_p \widetilde{Y}_j}{p} - \delta_{i,j}\frac{\mathrm{trace}(\Sigma_p)}{p}\right| > 3r(p)\right) \leq 2\exp(-K(\log(p))^{1+\varepsilon}).$$

And therefore,

$$P\left(\max_{i,j}\left|\frac{\widetilde{Y}_i' \Sigma_p \widetilde{Y}_j}{p} - \delta_{i,j}\frac{\mathrm{trace}(\Sigma_p)}{p}\right| > 3r(p)\right) \leq 2n^2\exp(-K(\log(p))^{1+\varepsilon}).$$

Using the Borel–Cantelli lemma, we reach the conclusion that

$$\max_{i,j}\left|\frac{\widetilde{Y}_i' \Sigma_p \widetilde{Y}_j}{p} - \delta_{i,j}\frac{\mathrm{trace}(\Sigma_p)}{p}\right| \leq 3r(p) = 3p^{2/m-1/2}\log(p) \qquad \text{a.s.}$$

Because the left-hand side is a.s. equal to $|\frac{Y_i' \Sigma_p Y_j}{p} - \delta_{i,j}\frac{\mathrm{trace}(\Sigma_p)}{p}|$, we reach the announced conclusion but with $r(p)$ replaced by $3r(p)$. Note that, of course, any multiple of $r(p)$, where the constant is independent of $p$, would work in the proof. In particular, by taking $\tilde{r}(p) = r(p)/3$, we reach the announced conclusion. $\square$

COROLLARY A.2. *Under the same assumptions as that of Lemma A.3, if we call $X_i = \Sigma_p^{1/2} Y_i$, we also have*

$$\max_{i \neq j}\left|\frac{\|X_i - X_j\|_2^2}{p} - 2\frac{\mathrm{trace}(\Sigma_p)}{p}\right| \leq p^{-1/2+2/m}(\log(p))^{(1+\varepsilon)/2} \qquad a.s.$$

PROOF. The proof follows immediately from the results of Lemma A.3, after we write

$$\|X_i - X_j\|_2^2 - 2\mathrm{trace}(\Sigma_p)$$
$$= [Y_i \Sigma_p Y_i - \mathrm{trace}(\Sigma_p)] + [Y_j \Sigma_p Y_j - \mathrm{trace}(\Sigma_p)] - 2Y_i' \Sigma_p Y_j.$$

Note that as explained in the proof of Lemma A.3, the constants in front of the bounding sequence do not matter, so we can replace $3p^{-1/2+2/m}(\log(p))^{(1+\varepsilon)/2}$ by $p^{-1/2+2/m}(\log(p))^{(1+\varepsilon)/2}$, and the result still holds. [In other words, we are really using Lemma A.3 with upper bound $p^{-1/2+2/m}(\log(p))^{(1+\varepsilon)/2}/3$.]
$\square$

LEMMA A.4. *Let $\{X_i\}_{i=1}^n$ be i.i.d. random vectors in $\mathbb{R}^p$ whose entries are i.i.d., mean 0, having the property that for 1-Lipschitz (with respect to Euclidean norm) functions $F$, if we denote by $m_F$ the median of $F(X_i)$,*

$$P(|F(X_i) - m_F| > r) \leq C\exp(-c(p)r^2),$$

where $C$ is independent of $p$ and $c$ is allowed to vary with $p$ (if it goes to zero, we assume it does so like $p^{-\alpha}$, $0 \leq \alpha < 1$). Call $\Sigma_p$ the covariance matrix of $X_1$. Assume that $\sigma_1(\Sigma_p)$ remains bounded in $p$. Then, under the triangular array construction of Theorem 2.3, we have, for any $\varepsilon > 0$,

$$\max_{i,j}\left|\frac{X_i'X_j}{p} - \delta_{i,j}\frac{\mathrm{trace}(\Sigma_p)}{p}\right| \leq (pc(p))^{-1/2}(\log(p))^{(1+\varepsilon)/2} \qquad a.s.$$

PROOF. The proof once again relies on concentration inequalities. First note that Proposition 1.11 combined with Proposition 1.7 in [29] shows that if $X_i$ and $X_j$ are independent and satisfy concentration inequalities with concentration function $\alpha(r)$ (with respect to Euclidean norm), then the vector $\binom{Y_i}{Y_j}$ also satisfies concentration inequalities with concentration function $2\alpha(r/2)$ with respect to Euclidean norm in $\mathbb{R}^{2p}$. (We note that Proposition 1.11 is proved for the metric on $\mathbb{R}^{2p}$ $\|\cdot\|_2 + \|\cdot\|_2$ where each Euclidean norm is a norm in $\mathbb{R}^p$, but the same proof goes through for Euclidean norm on $\mathbb{R}^{2p}$. Another argument would be to say that the metric $\|\cdot\|_2 + \|\cdot\|_2$ is equivalent to the norm of the full $\mathbb{R}^{2p}$ with the constants in the inequalities being 1 and $\sqrt{2}$ simply because for $a, b > 0$, $\sqrt{a^2 + b^2} \leq a + b \leq \sqrt{2}\sqrt{a^2 + b^2}$.)

Therefore, the arguments of Lemma A.2 go through without any problems with $\Sigma_p = \mathrm{Id}$ and $B_p^2 = 4/c(p)$. So a result similar to Corollary A.1 holds and we can apply the same ideas as in the proof of Lemma A.3 and get the announced result. □

COROLLARY A.3. *Under the assumptions of Lemma A.4, we have, for any $\varepsilon > 0$,*

$$\max_{i \neq j}\left|\frac{\|X_i - X_j\|_2^2}{p} - 2\frac{\mathrm{trace}(\Sigma_p)}{p}\right| \leq (pc(p))^{-1/2}(\log(p))^{(1+\varepsilon)/2} \qquad a.s.$$

PROOF. The proof is an immediate consequence of Lemma A.4, along the same lines as the proof of Corollary A.2. □

Finally, allow the same lines of proof; we have the following fact.

FACT A.1. *Let $\{X_i\}_{i=1}^n$ be i.i.d. random vectors in $\mathbb{R}^p$ whose entries are i.i.d., mean 0, having the property that for 1-Lipschitz (with respect to Euclidean norm) functions $F$, if we denote by $m_F$ the median of $F(X_i)$,*

$$P(|F(X_i) - m_F| > t) \leq C\exp(-c(p)t^b) \qquad \text{for some } b > 0,$$

*where $C$ is independent of $p$ and $c$ is allowed to vary with $p$ (if it goes to zero, we assume it does so like $p^{-\alpha}$, $0 \leq \alpha < b/2$). Call $\Sigma_p$ the covariance*

matrix of $X_1$. Assume that $\sigma_1(\Sigma_p)$ remains bounded in $p$. Then, we have, under the triangular array construction of Theorem 2.3, for any $\varepsilon > 0$,

$$\max_{i,j} \left| \frac{X_i' X_j}{p} - \delta_{i,j} \frac{\operatorname{trace}(\Sigma_p)}{p} \right| \leq (pc^{2/b}(p))^{-1/2} (\log(p))^{(1+\varepsilon)/b} \quad a.s.$$

Also, we then have

$$\max_{i \neq j} \left| \frac{\|X_i - X_j\|_2^2}{p} - 2\frac{\operatorname{trace}(\Sigma_p)}{p} \right| \leq (pc^{2/b}(p))^{-1/2} (\log(p))^{(1+\varepsilon)/b} \quad a.s.$$

The proof of this last fact follows the same step as that of Lemma A.4, with a slight adjustment since we need to replace 2 by $b$. For a related question and more details, we refer the reader to [19].

**(B) A linear algebraic result.** Finally, we finish this appendix with a linear algebraic lemma which we need in our approximations and is of independent interest.

LEMMA A.5. *Suppose M is real symmetric matrix with nonnegative entries. Suppose that E is a real symmetric matrix such that $\max_{i,j}|E_{i,j}| \leq \zeta$, for some $\zeta \geq 0$. Then, if $\sigma_1(A)$ is the largest singular value of matrix A and if $\circ$ represents the Hadamard product (i.e., entrywise multiplication of two matrices), we have*

$$\sigma_1(E \circ M) \leq \zeta \sigma_1(M).$$

PROOF. We first note that $E \circ M$ is real symmetric. Therefore,

$$\sigma_1(E \circ M) = \lim_{k \to \infty} [\operatorname{trace}((E \circ M)^{2k})]^{1/(2k)}.$$

Now we claim that

$$|\operatorname{trace}((E \circ M)^{2k})| \leq \zeta^{2k} \operatorname{trace}(M^{2k}).$$

To see this, recall that for a $p \times p$ matrix $A$,

$$\operatorname{trace}(A^k) = \sum_{1 \leq i_1, i_2, \ldots, i_k \leq p} A_{i_1, i_2} A_{i_2, i_3} \cdots A_{i_k, i_1}.$$

Now,

$$|A_{i_1,i_2} A_{i_2,i_3} \cdots A_{i_k,i_1}| \leq |A_{i_1,i_2}||A_{i_2,i_3}| \cdots |A_{i_k,i_1}|.$$

When $A = E \circ M$, $A_{i,j} = E_{i,j} M_{i,j}$. Since $M_{i,j} \geq 0$, we therefore have $|E_{i,j} \times M_{i,j}| \leq \zeta M_{i,j}$. Hence,

$$|\operatorname{trace}((E \circ M)^k)| \leq \sum_{1 \leq i_1, i_2, \ldots, i_k \leq p} \zeta^k M_{i_1,i_2} M_{i_2,i_3} \cdots M_{i_k,i_1} = \zeta^k \operatorname{trace}(M^k).$$

So

$$[\text{trace}((E \circ M)^{2k})]^{1/(2k)} \leq \zeta [\text{trace}(M^{2k})]^{1/(2k)}.$$

Taking limits as $k \to \infty$ concludes the proof. $\square$

**Acknowledgments.** I would like to thank Bin Yu for stimulating my interest in the questions considered in this paper and for interesting discussions on the topic. I would like to thank Elizabeth Purdom for discussions about kernel analysis and Peter Bickel for many stimulating discussions about random matrices and their relevance in statistics. I would also like to thank an anonymous referee for useful and constructive comments that resulted in an improved presentation of the paper.

DEPARTMENT OF STATISTICS
UNIVERSITY OF CALIFORNIA, BERKELEY
367 EVANS HALL
BERKELEY, CALIFORNIA 94720-3860
USA
E-MAIL: nkaroui@stat.berkeley.edu